\documentclass[12pt,a4paper]{amsart}
\usepackage[cm,headings]{fullpage}
\usepackage{amssymb,amscd}

\theoremstyle{plain}
\newtheorem{theorem}{Theorem}[section]
\newtheorem{proposition}[theorem]{Proposition}
\newtheorem{lemma}[theorem]{Lemma}
\newtheorem{corollary}[theorem]{Corollary}

\theoremstyle{remark}
\newtheorem*{remark}{Remark}
\newtheorem*{remarks}{Remarks}

\def\qed{$\blacksquare$ \bigskip}

\newcommand{\set}[2]{\ensuremath{\{ #1 \>|\> #2 \}}}

\def\liebrack{\ensuremath{[\,\cdot\, , \cdot\,]}}      
\def\curlybrack{\ensuremath{\{\,\cdot \,, \cdot\,\}}}  

\begin{document}

\title[Deformations of $W_1(n) \otimes A$]
{Deformations of $W_1(n) \otimes A$ and modular semisimple Lie algebras with
a solvable maximal subalgebra}
\author{Pasha Zusmanovich}
\address{}
\email{justpasha@gmail.com}
\date{March 30, 2002; last minor revision September 16, 2014}
\thanks{\textsf{arXiv:math/0204004}; J. Algebra \textbf{268} (2003), 603--635}


\begin{abstract}
In one of his last papers, Boris Weisfeiler proved that if a modular semi\-simple Lie algebra 
possesses a solvable maximal subalgebra which defines in it a long filtration, 
then the associated graded algebra is isomorphic to the one constructed from the 
Zassenhaus algebra tensored with the divided powers algebra.
We completely determine such class of algebras, calculating in the process 
low-dimensional cohomology groups of the Zassenhaus algebra tensored with any 
associative commutative algebra.
\end{abstract}

\maketitle

\section*{Introduction}

The ultimate goal of this paper is to describe semisimple finite-dimensional Lie 
algebras over an algebraically closed field of characteristic $p > 5$, having a solvable 
maximal subalgebra (that is, a maximal subalgebra which is solvable)
which determines a ``long filtration'', as defined below. 
We hope, however, that some intermediate results contained here are of independent 
interest.

Simple Lie algebras with a solvable maximal subalgebra were described by 
B. Weisfeiler \cite{W}. We depend heavily on his results and to a certain extent 
this paper may be considered as a continuation of the Weisfeiler's paper. 
Since its appearance the classification of modular simple Lie algebras has been 
completed (announced in \cite{SW} and elaborated in a series of papers among which
\cite{Str} is the last one), and recently the approach to the classification problem 
has been reworked in a series of papers among which \cite{PS} is the latest,
including the low characteristic cases.

Though the classification provides a powerful tool for solutions of many problems in 
modular Lie algebras theory, the question considered here remains non-trivial even 
modulo this classification. Moreover, we hope that the result we obtain here 
in particular, and the cohomological technique we use to prove it in general, 
may simplify to certain degree the classification itself.

Let us recall the contents of Weisfeiler's paper. 
He considers semisimple modular Lie algebra $\mathfrak L$ with a solvable 
maximal subalgebra 
$\mathfrak L_0$. $\mathfrak L_0$ defines a filtration in $\mathfrak L$ via 
$\mathfrak L_{i+1} = \set {x \in \mathfrak L_i}{[ x, \mathfrak L ] = \mathfrak L_i}$
(though in general the filtration can be prolonged also to the negative side, 
in the case under consideration we can let $\mathfrak L_{-1} = \mathfrak L$).

When the term $\mathfrak L_1$ of this filtration does not vanish, 
the filtration is called {\it long}, otherwise it is called {\it short}. 
Weisfeiler proved that when the filtration is long, the associated graded algebra 
is isomorphic to $S \otimes O_m + 1 \otimes \mathfrak D$, where $S$ coincides either 
with $sl(2)$, the three-dimensional simple algebra, or with $W_1(n)$, the
Zassenhaus algebra, $O_m$ is the reduced polynomial ring in $m$ variables, 
and $\mathfrak D$ is a derivation algebra of $O_m$. The grading is ``thick'' 
in the sense that it is completely determined by the standard grading of 
$W_1(n)$ or $sl(2)$, containing therefore the whole tensor factor $O_m$ in each 
component. 
In the short filtration case, Weisfeiler proved that the initial algebra
$\mathfrak L$ possesses a $\mathbb Z_p$-grading with very restrictive conditions. 
Then, considering the case of simple $\mathfrak L$, he derived that in the long filtration
case, $\mathfrak L$ is isomorphic either to $sl(2)$ or to $W_1(n)$
(in fact, this follows immediately from the results of Kuznetsov \cite{K}, which are also 
important for us here), and the short filtration case does not occur.

Here we study the long filtration case. We determine all filtered algebras whose 
associated graded algebra is $W_1(n) \otimes O_m + 1 \otimes \mathfrak D$ or 
$sl(2) \otimes O_m + 1 \otimes \mathfrak D$ with the above-mentioned ``thick'' grading. 
This is done in the framework of the deformation theory due to Gerstenhaber.
In this theory the second cohomology group of a Lie algebra with coefficients in
the adjoint module plays a significant role (for an excellent account of this subject,
see [GS]). As it turns out that the ``tail'' $1 \otimes \mathfrak D$ is not important 
in these considerations, one needs to compute 
$H^2 (W_1(n) \otimes O_m, W_1(n) \otimes O_m)$. 
In a (slightly) more general setting, we compute 
$H^2 (W_1(n) \otimes A, W_1(n) \otimes A)$ for an arbitrary associative commutative 
algebra $A$ with unit. ($H^2 (sl(2) \otimes A, sl(2) \otimes A)$) was earlier
computed by Cathelineau [C]). This calculation seems sufficiently interesting for
its own sake, as a nontrivial example of the low-dimensional cohomology of current
Lie algebras $L \otimes A$. This also may be considered as a complement to the 
Cathelineau's computation of the second cohomology group of the current Lie algebra
$\mathfrak g \otimes A$ extended over a classical simple Lie algebra $\mathfrak g$, 
as well as generalization of Dzhumadil'daev-Kostrikin computations of 
$H^2 (W_1(n), W_1(n))$ [DK].

The knowledge of $H^2 (W_1(n) \otimes O_m, W_1(n) \otimes O_m)$ allows to solve 
the problem of determining all filtered algebras associated with the graded structure 
mentioned above. The answer is not very surprising -- all such algebras have a socle 
isomorphic to $W_1(k) \otimes O_l$ for some $k$ and $l$.

The contents of this paper are as follows. \S 1 contains some preliminary material, 
the most significant of which is a representation of $W_1(n)$ as a deformation of
$W_1(1) \otimes O_{n-1}$, due to Kuznetsov \cite{K}. It turns out that it is much easier 
to perform cohomological calculations using this representation. 
Following Kuznetsov, we define a class of Lie algebras $\mathfrak L(A, D)$ 
which are certain deformations of $W_1(1) \otimes A$ defined by means of a derivation
$D$ of $A$. Then we compute $H^2 (\mathfrak L(A, D), \mathfrak L(A, D))$
in two steps: first, in \S 2, we compute $H^2 (W_1(1) \otimes A, W_1(1) \otimes A)$, 
and then, in \S 3, we determine $H^2 (\mathfrak L(A, D), \mathfrak L(A, D))$, 
using a spectral sequence abutting to $H^* (\mathfrak L(A, D), \mathfrak L(A, D))$ 
with the $E_1$-term isomorphic to $H^*(W_1(1) \otimes A, W_1(1) \otimes A)$.

Parallel to the results for the second cohomology group, we state similar results
for the first cohomology group, as well as for the second cohomology group with
trivial coefficients, which are useful later, in \S 5.

In \S 4, using the Kuznetsov's isomorphism, we transform the results about \allowbreak
${H^2 (\mathfrak L(A,D), \mathfrak L(A,D))}$ into those about 
$H^2 (W_1(n) \otimes A, W_1(n) \otimes A)$. This section contains also all necessary 
computations related to reduced polynomial rings, particularly, of their
Harrison cohomology. After that, in \S 5 we formulate a theorem 
about filtered deformations of $W_1(n) \otimes A + 1 \otimes \mathfrak D$ and of 
$sl(2) \otimes A + 1 \otimes \mathfrak D$ and derive it almost immediately from 
preceding results. It turns out that each such deformation strictly related to 
the class $\mathfrak L(A, D)$ (for a different $A$), so in \S 6 we determine
all semisimple algebras in this class up to isomorphism, completing therefore 
the consideration of the long filtration case (Theorem 6.4).

Since the present paper is overloaded with different kinds of computations, we 
omit some of them which are similar to those already presented, or just too tedious. 
We believe that this will not cause inconvenience to the reader.

\section{Preliminaries}

In this section we recall all necessary notions, notation, definitions, results and
theories, as well as define a class of algebras $\mathfrak L(A, D)$ important for 
further considerations.

The ground field $K$ is assumed to be of characteristic $p > 3$, unless
otherwise is stated explicitly.
(When appealing to Weisfeiler's results, we have to assume the ground
field is algebraically closed of characteristic $p > 5$).

As we deal with modular Lie algebras, it is not surprising that the divided powers
algebra $O_1(n)$ plays a significant role in our considerations. 
Recall that $O_1(n)$ is the commutative associative algebra with basis
$\set {x^i}{0 \le i < p^n}$ and multiplication $x^ix^j = \binom {i+j}j x^{i+j}$. 
It is isomorphic to the reduced polynomial ring 
$O_n = K[x_1, \dots, x_n]/(x_1^p, \dots , x_n^p)$, the isomorphism is given by
\begin{equation}\tag{1.1}
x_1^{\alpha_1} x_2^{\alpha_2} \dots x_n^{\alpha_n} \mapsto 
\alpha_1! \alpha_2! \dots \alpha_n!
x^{\alpha_1 + p\alpha_2 + \dots + p^{n-1}\alpha_n} .
\end{equation}

The subalgebra $\set {x^i}{1 \le i < p^n}$, denoted as $O_1(n)^+$ (or $O_n^+$)
is the single maximal ideal of $O_1(n)$. The invertible elements of $O_1(n)$ are exactly those 
not lying in $O_1(n)^+$.

The Zassenhaus algebra $W_1(n)$ is the Lie algebra of derivations of $O_1(n)$ of the kind
$u\partial$, where $u \in O_1(n)$ and $\partial(x^j) = x^{j-1}$. It possesses a basis 
$\set {e_i = x^{i+1}\partial}{-1 \le i \le p^n-2}$ with bracket
\begin{equation}\notag
[e_i, e_j] = 
N_{ij}e_{i+j}, \qquad\text{where }N_{ij} = \binom {i + j + 1}j - \binom {i + j + 1}i .
\end{equation}

The grading $W_1(n) = \bigoplus_{i=-1}^{p^n-2} Ke_i$ is called \textit{standard}. 
In the case $n = 1$ it coincides with the root space decomposition relative to the action of
the semisimple element $e_0$.

Notice the following properties of the coefficients $N_{ij}$:

\begin{list}{}{}
\item[(1)] $N_{ij} = 0 \text{ if } -1 \le i,j \le p-2, i+j \ge p-2$
\item[(2)] $N_{ij} = N_{i-1,j} + N_{i,j-1}$
\item[(3)] $N_{i,j-p} = N_{ij} \text{ if } -1 \le i \le p-2 \le j$
\end{list}

The first two are obvious, the third may be found, for example, in \cite{DK}. Notice
also that if $0 \le i \le j \le p$, then $\binom{p-i}{p-j} = (-1)^{j-i}\binom{j-1}{i-1}$
and if $i = \sum_{n \ge 0} i_np^n$, $j = \sum_{n \ge 0} j_np^n$ are $p$-adic decompositions, 
then
\begin{equation}\notag
\binom ij = \prod_{n \ge 0} \binom{i_n}{j_n}
\end{equation}
(the latter is known as Lucas' theorem).

The derivation algebra $Der(W_1(n))$ is generated (linearly) by inner derivations
of $W_1(n)$ and derivations $(ad \, e_{-1})^{p^t}, 1 \le t \le n-1$ 
(so the latters constitute a basis of $H^1 (W_1(n), W_1(n))$
(cf. \cite{B1} or \cite{D2}).

The whole derivation algebra of $O_1(n) \simeq O_n$, known as a general Lie algebra of
Cartan type $W_n$, is freely generated as $O_1(n)$-module by 
$\set{ \partial^{p^i}}{0 \le i \le n-1}$ (or, in terms of $O_n$, by 
$\set{ \partial / \partial x_i}{1 \le i \le n}$)
(cf. \cite{BIO}, \cite{K}, or \cite{W}).

Let $L$ be a Lie algebra and $A$ an associative commutative algebra with unit. The
Lie structure on the tensor product $L \otimes A$ is defined via 
$[x \otimes a, y \otimes b] = [x, y] \otimes ab$. If $\mathfrak D$ is a subalgebra of 
$Der(A)$, then $L \otimes A + 1 \otimes \mathfrak D$ is defined as a semidirect product
where $1 \otimes \mathfrak D$ acts on $L \otimes A$ by 
$[x \otimes a, 1 \otimes d] = x \otimes d(a)$. For $a \in A$, $R_a$ stands for the
multiplication on $a$ in $A$.

We will need the following elementary results.

\begin{proposition}\label{1.1}
\hfill
\begin{enumerate}
\item $Z (L \otimes A) = Z(L) \otimes A$
\item $(1 \otimes Der(A)) \cap ad(L \otimes A) = 0$.
\end{enumerate}
\end{proposition}

\begin{proof}
(i) Obviously $Z(L) \otimes A \subseteq Z(L \otimes A)$. 
Let $\sum z_i \otimes a_i \in Z(L \otimes A)$. We may assume that all $a_i$ are 
linearly independent. Then
\begin{equation}\notag
[\, \sum z_i \otimes a_i, x \otimes 1 \,] = \sum [z_i,x] \otimes a_i = 0
\end{equation}
for every $x \in L$, which together with our assumption implies $z_i \in Z(L)$ for all $i$.

(ii) Let $1 \otimes d = \sum ad\,x_i \otimes a_i \in (1 \otimes Der(A)) \cap ad(L \otimes A)$. 
Applying it to $y \otimes 1$, we get $\sum [y, x_i] \otimes a_i = 0$ whence $x_i \in Z(L)$ 
for all $i$, and $d = 0$. 
\end{proof}

\noindent \textit{Definition}.
Let $D \in Der(A)$. Define $\mathfrak L(A, D)$ 
to be a Lie algebra with the underlying vector space $W_1(1) \otimes A$ and Lie bracket
$\{x,y\} = [x,y] + \Phi_D (x,y)$, where $\liebrack$ is the ordinary bracket on 
$W_1(1) \otimes A$, and
\begin{equation}\tag{1.2}
\Phi_D (e_i \otimes a, e_j \otimes b) = \begin{cases} 
e_{p-2} \otimes (aD(b) - bD(a)), & i = j = -1      \\
0,                               & \text{otherwise}                        
\end{cases}
\end{equation}
for $a,b \in A$.
\bigskip

\begin{remark}[Referee] 
In \cite{ree}, Ree considered a class of Lie algebras
which are subalgebras of $Der(B)$ for a commutative associative algebra $B$
(freely) generated over $B$. The algebras $\mathfrak L(A, D)$ belong to this class. 
Indeed, let $D\in Der(A)$ and consider
$d = \partial\otimes 1 + x^{p-1}\otimes D \in Der(O_1\otimes A)$. One can easily see,
by identifying $O_1$ and $W_1(1)$ as vector spaces via $x^i \mapsto e_{i-1}$, $0\le i < p$,
that $\mathfrak L(A, D)$ is nothing else than the Lie algebra of derivations of $O_1\otimes A$
of the form $\set{bd}{b \in O_1\otimes A}$ (i.e. freely generated,
as a module over $O_1\otimes A$, by a single derivation $d$).
\end{remark}

The following is crucial for our considerations.

\begin{proposition}[Kuznetsov]\label{kuznetsov}
$W_1(n) \otimes A \simeq \mathfrak L(O_1(n -1) \otimes A, \partial \otimes 1)$.
\end{proposition}

\begin{proof}
This obviously follows from the isomorphism 
$W_1(n) \simeq \mathfrak L(O_1(n-1), \partial)$, noticed in \cite{K}.
A direct calculation shows that the mapping
\begin{equation}\notag
e_{pk+i} \mapsto e_i \otimes x^k, \qquad -1 \le i \le p-2, \quad 0 \le k < p^{n-1}
\end{equation}
provides the isomorphism desired.
\end{proof}

The reason why we prefer to deal with such realization of $W_1(n)$ lies in the fact
that $e_0 \otimes 1$ remains a semisimple element in $\mathfrak L(A, D)$ with root spaces 
$e_i \otimes A$. So we obtain the grading of length $p$, and not of length $p^n$ 
as in the case of $W_1(n) \otimes A$.
The significance of the "good"(=short) root space decomposition follows from the
well-known theorem about the invariance of the Lie algebra cohomology under the
torus action.

Introduce a filtration
\begin{equation}\tag{1.3}
\mathfrak L(A, D) = \mathfrak L_{-1} \supset \mathfrak L_0 \supset \mathfrak L_1 \supset
\dots \supset \mathfrak L_{p-2}
\end{equation}
by $\mathfrak L_i = \bigoplus_{j \ge i} e_j \otimes A$.

In general, for a given decreasing filtration $\{\mathfrak L_i\}$ of a Lie algebra 
$\mathfrak L$, $gr \mathfrak L = \bigoplus_i \mathfrak L_i / \mathfrak L_{i+1}$ 
will denote its associated graded algebra.

The following is evident.

\begin{proposition}
The graded Lie algebra $gr \mathfrak L(A, D)$ associated with filtration
(1.3), is isomorphic to $W_1(1) \otimes A = \bigoplus_{i=-1}^{p-2} e_i \otimes A$.
\end{proposition}

This is the place where deformation theory enters the game. It is known that each
filtered algebra can be considered as a deformation of its associated graded algebra
$L = \bigoplus L_i$ (for this fact as well as for all necessary background in the deformation
theory we refer to \cite{GS}). One calls such deformations \textit{filtered deformations}
(or $\{ L_i \}$-deformations in the terminology of \cite{DK}). 
As the space of infinitesimal deformations coincides with the cohomology group 
$H^2(L, L)$, the space of infinitesimal filtered deformations coincides with its subgroup 
$H_+^2 (L, L) = 
\set{ \overline\phi \in H^2(L, L)}{ \phi (L_i, L_j) \subset \bigoplus_{k \ge 1} L_{i+j+k}}$. 
To describe all filtered deformations, one needs to investigate prolongations of 
infinitesimal ones, obstructions to which are described by Massey
products $[\phi, \psi] \in H^3(L,L)$ defined as
\begin{equation}\notag
[\phi, \psi](x,y,z) = \phi(\psi(x,y),z) + \psi(\phi(x,y),z) + \curvearrowright . 
\end{equation}
(this product arises from the graded Lie (super)algebra
structure on $H^*(L, L)$). 

We formulate just a small part of this broad subject needed for our purposes.

\begin{proposition}[cf. \cite{GS} or direct verification]
Let $L$ be a finitely graded Lie algebra such that Massey product of any two elements 
of $Z_+^2(L, L)$ is zero. 
Then any filtered Lie algebra $\mathfrak L$ such that $gr\,\mathfrak L \simeq L$
(as graded algebras), is isomorphic to a Lie algebra with underlying vector space
$L$ and Lie bracket $\curlybrack = \liebrack + \Phi$ for some $\Phi \in Z_+^2(L,L)$.
\end{proposition}

Note in that connection that $[\Phi_D, \Phi_D] = 0$. We will see later that this holds also
for other ``positive'' 2-cocycles on $W_1(1)\otimes A$ (and more generally, on $W_1(n)\otimes A$),
so Proposition 1.4 will be applicable in our situation.

Now we formulate the Weisfeiler's main result \cite{W}:

\begin{theorem}[Weisfeiler]
(The ground field $K$ is algebraically closed of characteristic $p > 5$).

Let $\mathfrak L$ be a semisimple Lie algebra with a solvable maximal subalgebra
$\mathfrak L_0$. Suppose that $\mathfrak L_0$ defines a long filtration 
in $\mathfrak L$.
Then $\mathfrak L$ is a filtered deformation of a graded Lie algebra 
$L = S \otimes O_m + 1 \otimes \mathfrak D$, where $S = sl(2)$ or $W_1(n)$ equipped
with the standard grading $\bigoplus_i Ke_i$, $\mathfrak D \subset Der(O_m)$,
and the graded components are:
\begin{equation}\notag
L_i = \begin{cases}
e_0 \otimes O_m + 1 \otimes \mathfrak D, & i=0 \\
e_i \otimes O_m,                         & i \ne 0.
\end{cases}
\end{equation}
\end{theorem}

Further, the Harrison cohomology $Har^* (A, A)$ with coefficients in the adjoint
module $A$ plays a role in our considerations. Note that $Har^1 (A, A) = Der(A)$
and Harrison 2-cocycles, denoted by $\mathcal Z^2(A, A)$, are just symmetrized Hochschild
2-cocycles (cf. \cite{harrison} where this cohomology was introduced and \cite{GS} for a more modern
treatment). $\delta$ refers to the Harrison (=Hochschild) coboundary operator, i.e.
\begin{align}
\delta G(a, b)  &= aG(b) + bG(a) - G(ab)                 \notag   \\
\delta F(a,b,c) &= aF(b,c) - F(ab,c) + F(a,bc) - F(a,b)c \notag
\end{align}
for $G \in Hom(A, A)$ and $F \in Hom(A \otimes A, A)$. The action of $Der(A)$ on $Har^2(A, A)$
is defined via
\begin{equation}\notag
D \star F(a, b) = F(D(a), b) + F(a, D(b)) - D(F(a, b)).
\end{equation}
The same formula defines the action of $Der(L)$ on the cohomology $H^2(L, L)$ of the
Lie algebra $L$.

Considering the $L$-action on $H^*(L, L)$, the well-known fact says that if $T$ is an
abelian subalgebra relative to which $L$ decomposes into a sum of eigenspaces 
$L = \bigoplus L_{\alpha}$, then one can decompose the complex into the sum of subcomplexes
\begin{equation}\notag
C_{\alpha}^n = 
\set{\phi \in C^n(L, L)}{\phi(L_{\alpha_1}, \dots, L_{\alpha_n}) \subseteq 
                              L_{\alpha_1 + \dots + \alpha_n + \alpha}}
\end{equation}
and, moreover, $H^*(C_{\alpha}) = 0$ for $\alpha \ne 0$ (cf. \cite{F}, Theorem 1.5.2).

Similarly, any $\mathbb Z$-grading $L = \bigoplus L_i$ induces a $\mathbb Z$-grading 
on the cohomology group $H^*(L, L)$, as the initial complex $C^*(L, L)$ splits 
into the sum of subcomplexes $C_i^*(L, L)$, where
\begin{equation}\tag{1.4}
C_i^n(L,L) = 
\set{\phi \in C^n(L,L)}{\phi(L_{i_1}, \dots, L_{i_n}) \subseteq L_{i_1 + \dots + i_n +i}}.
\end{equation}
The cocycles, coboundaries and cohomology of these subcomplexes form the modules 
denoted by $Z_i^n(L, L)$, $B_i^n(L, L)$ and $H_i^n(L, L)$ respectively. 
If there is an element $e \in L$ whose action on $L_i$ is multiplication by $i$
then $H^*(C_i) = 0$ for $i \ne 0 \,mod\,p$.

The symbol $\curvearrowright$ after an expression refers to the sum of all cyclic permutations
(in $S(3)$) of letters and indices occuring in that expression.

\section{Low-dimensional cohomology of $W_1(1) \otimes A$}

The aim of this section is to establish the following isomorphisms.

\begin{proposition}\label{2.1}
\hfill
\begin{enumerate}
\item $H^1(W_1(1) \otimes A, W_1(1) \otimes A) \simeq Der(A)$
\item $H^2(W_1(1) \otimes A, W_1(1) \otimes A)$
\flushright{
$\simeq H^2(W_1(1), W_1(1)) \otimes A \>\oplus\> Der(A) \>\oplus\> Der(A) 
\>\oplus\> Har^2(A,A)$.
}
\end{enumerate}
\end{proposition}

Before beginning the proof, let us make several remarks.

Part (i) follows from \cite{block-diff}, Theorem 7.1 (formulated in terms
of derivation algebras). Alternatively, one may prove it in a similar (and much
easier) way as (ii). 
Perhaps it should be remarked only that the basic 1-cocycles on $W_1(1) \otimes A$ 
can be given as $1 \otimes D$ for $D \in Der(A)$.

So we will concentrate our attention on (ii).

The cohomology group $H^2 (W_1(n), W_1(n))$ was computed in \cite{DK}. Particularly,
\newline
$\dim H^2 (W_1(1), W_1(1)) = 1$ and the single basic cocycle can be chosen as:
\begin{equation}\tag{2.1}
\phi (e_i, e_j) = \begin{cases}
N_{ij}/p \cdot e_{i+j-p}, & i+j \ge p-1 \\
0,                        & \text{otherwise}
\end{cases}
\end{equation}
where $N_{ij}/p$ denotes a (well defined) element of the field $K$ which is obtained from
$N_{ij}$ by division by $p$ and further reduction modulo $p$.

The appearance of the first and last terms in (ii) is evident: the corresponding
parts of the cohomology group are spanned by the classes of cocycles
\begin{align}
\Theta_{\phi, u}: x \otimes a \wedge y \otimes b &\mapsto \phi(x, y) \otimes abu \tag{2.2} \\
\Upsilon_F:       x \otimes a \wedge y \otimes b &\mapsto [x,y] \otimes F(a,b)   \tag{2.3}
\end{align}
respectively, where $\phi \in Hom(L \otimes L, L)$, $u \in A$, and 
$F \in Hom(A \otimes A, A)$.
We will denote the cochains of type (2.2) with $u = 1$ as $\Theta_{\phi}$ 
(so actually $\Theta_{\phi, u} = (1 \otimes R_u) \circ \Theta_{\phi}$).

We have the following simple proposition.

\begin{proposition}\label{2.2}
Let $L$ be a Lie algebra which is not 2-step nilpotent. Then
\begin{enumerate}
\item 
$(1 \otimes R_u) \circ \Theta_{\phi} \in Z^2(L \otimes A, L \otimes A)$ 
if and only if either $\phi \in Z^2(L, L)$ or $u=0$

\item $\Upsilon_F \in Z^2(L \otimes A, L \otimes A)$ 
if and only if $F \in \mathcal Z^2(A, A)$.
\end{enumerate}
\end{proposition}

\begin{proof}
We will prove the second part only, the first one is similar. The cocycle
equation for $\Upsilon_F$ together with Jacobi identity gives
\begin{equation}\tag{2.4}
[[x,y],z] \otimes \delta F(a,c,b) + [[z,x],y] \otimes \delta F(a,b,c) = 0.
\end{equation}
Since $[[L,L],L] \ne 0$ and $p\ne 3$, one may choose $x, y \in L$ 
such that $[[y, x], x] \ne 0$. Setting $z = x$, one gets 
$F \in \mathcal Z^2(A, A)$. Conversely, the last condition implies (2.4).
\end{proof}

It is possible to prove also that for any Lie algebra $L$ these cocycles are 
cohomologically independent, whence $H^2(L \otimes A, L \otimes A)$ must contain 
$H^2(L, L) \otimes A$ and $Har^2(A, A)$ as direct summands.

Let us define now explicitly the remaining classes of basic cocycles: 
$\Phi_D$ is already defined by (1.2), and
\begin{multline}\tag{2.5}
\Psi_D (e_i \otimes a, e_j \otimes b) \\ = \begin{cases}
e_{i+j} \otimes (\binom{i+j+1}j bD(a) - \binom{i+j+1}i aD(b)), & -2 < i+j < p-1 \\
0                                                              & \text{otherwise}.
\end{cases}
\end{multline}

\begin{lemma}\label{2.3}
For any $D \in Der(A)$, $\Psi_D, \Phi_D \in Z^2(W_1(1) \otimes A, W_1(1) \otimes A)$.
\end{lemma}

\begin{proof}
We perform necessary calculations for $\Psi_D$, leaving the easier case of 
$\Phi_D$ to the reader 
(in fact, that $\Phi_D$ is a 2-cocycle on $W_1(1) \otimes A$ follows from the Jacobi identity
in $\mathfrak L(A, D)$).

Isolating the coefficient of $e_{i+j+k}\otimes abD(c)$ in the cocycle equation
for $\Psi_D$, we get
\begin{multline}\notag
- N_{ij}    \binom{i+j+k+1}{i+j}
+ N_{jk}    \binom{i+j+k+1}i     
+ N_{ki}    \binom{i+j+k+1}j          \\
+ N_{i,j+k} \binom{j+k+1}j
- N_{j,k+i} \binom{k+i+1}i        = 0.
\end{multline}
The last relation can be verified immediately.
\end{proof}

The element $e_0 \otimes 1$ acts semisimply on $W_1(1) \otimes A$, as well as on 
$\mathfrak L(A, D)$. The roots of $ad(e_0 \otimes 1)$-action lie in the prime subfield
and the root spaces are:
\begin{equation}\tag{2.6}
L_{[i]} = e_i \otimes A, \quad [i] \in \mathbb Z_p, \quad -1 \le i \le p-2.
\end{equation}
Thus any cocycle in $Z^2 (W_1(1) \otimes A, W_1(1) \otimes A)$ is cohomologous to 
a cocycle in
\begin{equation}\tag{2.7}
Z_{[0]}^2 (W_1(1) \otimes A, W_1(1) \otimes A) = Z_{-p}^2 \oplus Z_0^2 \oplus Z_p^2
\end{equation}
as noted above.

\begin{lemma}\label{2.4}
Let $\{u_i\}$ be linearly independent elements of A, $\{D_i\}$ be linearly
independent derivations of $A$, $\{F_i\}$ be cohomologically independent cocycles in
$\mathcal Z^2(A, A)$.

Then cocycles 
$(1 \otimes R_{u_i}) \circ \Theta_{\phi}$, 
$\Psi_{D_i}$,
$\Upsilon_{F_i}$,
$\Phi_{D_i}$
(defined in (2.2), (2.5), (2.3) and (1.2) respectively), are cohomologically independent.
\end{lemma}

\begin{proof}
As the cocycles of the first type belong to the ($-p$)th component of 
$Z^2 (W_1(1) \\ \otimes A, W_1(1) \otimes A)$, 
the cocycles of the second and third type --
to the zero component, the cocycles of the fourth type -- to the $p$th component, 
and the degree of any coboundary is in the range between $1-p$ and $p-1$,
one needs only to show the independence of cocycles of the form $\Psi_{D_i}$ and
$\Upsilon_{F_i}$.

Suppose that there is a linear combination of the above-mentioned cocycles equal
to a coboundary $d\omega$. Clearly this condition can be written as
\begin{equation}\tag{2.8}
\Psi_D + \Upsilon_F = d\omega
\end{equation}
where $D$, $F$ are some linear combinations of $D_i$'s and $F_i$'s, respectively.

Due to the $e_0 \otimes 1$-action on $\mathfrak L(A, D)$, we may assume that $\omega$
preserves the root space decomposition (2.6), i.e.
\begin{equation}\notag
\omega (e_i \otimes a) = e_i \otimes X_i(a)
\end{equation}
for some $X_i \in Hom(A, A)$.

Evaluating the left and right sides of (2.8) for the pair $e_0 \otimes a$, $e_0 \otimes 1$, 
one gets $D = 0$. Then (2.8) reduces to
\begin{equation}\tag{2.9}
F(a, b) = aX_j(b) + bX_i(a) - X_{i+j}(ab)
\end{equation}
for all $i,j$ such that $N_{ij} \ne 0$.

Substituting in (2.9) $j = 0$ and using the symmetry of $F$, we get $F = \delta X_0$. Since
$F$ is a linear combination of cohomologically independent Harrison cocycles, $F = 0$.
We see that all elements entering (2.8) vanish, whence all coefficients in the initial
linear combinations of cocycles are equal to zero.
\end{proof}

Now, to prove Proposition 2.1(ii), one merely needs to show that each cocycle
$\phi \in Z_{[0]}^2 (W_1(1) \otimes A, W_1(1) \otimes A)$ is cohomologous to the sum of 
the previous cocycles.

Let
\begin{equation}\tag{2.10}
\phi = \phi_{-p} + \phi_0 + \phi_p, \quad 
\phi_k (e_i \otimes A, e_j \otimes A) \subseteq e_{i+j+k} \otimes A, \quad k=-p, 0, p
\end{equation}
be a decomposition corresponding to (2.7). It is immediate that 
$d\phi = 0 \iff d\phi_{-p} = d\phi_0 = d\phi_p = 0$.

The next three lemmas elucidate the form of cocycles $\phi_{-p}, \phi_0, \phi_p$ respectively.
Two of them are formulated in a slightly more general setting which will be used
later, in \S 3.

\begin{lemma}\label{2.5}
$\phi_{-p} = (1 \otimes R_u) \circ \Theta_{\phi}$ for some $u \in A$.
\end{lemma}

\begin{proof}
Write
\begin{equation}\notag
\phi_{-p} (e_i \otimes a, e_j \otimes b) = \begin{cases}
e_{i+j-p} \otimes X_{ij}(a,b), & i+j \ge p-1       \\
0,                             & \text{otherwise}
\end{cases}
\end{equation}
for certain $X_{ij} \in Hom(A \otimes A, A)$. Writing the cocycle equation for triples
$e_i \otimes a, e_j \otimes b, e_{-1} \otimes 1$ and 
$e_i \otimes a, e_j \otimes 1, e_0 \otimes c$, one obtains respectively:
\begin{align}
X_{ij}(a,b) &= X_{i-1,j}(a,b) + X_{i,j-1}(a,b), \quad i+j>p-1                      \tag{2.11} \\
X_{ij}(a,c)  &= \frac{i+j}j cX_{ij}(a,1) - \frac ij X_{ij}(ac,1), \quad i+j \ge p-1 \tag{2.12}
\end{align}
The last equality in the case $i+j = p$ entails
\begin{equation}\notag
X_{ij}(a,c) = X_{ij}(ac,1), \quad i+j=p.
\end{equation}
Now the cocycle equation for the triple 
$e_i \otimes a, e_j \otimes b, e_1 \otimes 1, i+j = p-1, i,j \ne 1$ implies
\begin{equation}\notag
X_{ij}(a,b) = - N_{1,i}(X_{i+1,j}(a, b) + X_{i,j+1}(a, b)), \quad i+j = p-1.
\end{equation}
Substitution of the last but one equality into the last one yields
\begin{equation}\tag{2.13}
X_{ij}(a, b) = Y_{ij}(ab), \quad i+j = p-1, \quad i,j \ne 1
\end{equation}
where $Y_{ij} = - N_{1,i}(X_{i+1,j}(a, 1) + X_{i,j+1}(a, 1))$. 
Substituting this in its turn, in (2.12) (with $i + j = p-1$), one gets
\begin{equation}\notag
Y_{ij}(ac) = cY_{ij}(a)
\end{equation}
which implies $Y_{ij}(a) = au_{ij}$ for some $u_{ij} \in A$. Hence
\begin{equation}\notag
X_{ij}(a,b) = abu_{ij}, \quad i+j = p-1, \quad i,j \ne 1.
\end{equation}
Writing the cocycle equation for the triple $e_1 \otimes a, e_1 \otimes 1, e_{p-2} \otimes 1$, 
one obtains
\begin{equation}\notag
X_{1,p-2}(a, 1) = aX_{1,p-2}(1, 1).
\end{equation}
Substituting this in (2.12) under the particular case $i=1, j=p-2$, one deduces (2.13)
also in this case, with $u_{1,p-2} = X_{1,p-2}(1, 1)$. Then writing the cocycle equation for
triple $e_i \otimes 1, e_j \otimes 1, e_1 \otimes 1, i+j = p-2, i,j \ne 0$, 
and taking into account (2.13), one obtains
\begin{equation}\notag
N_{1,i}u_{i+1,j} + N_{1,j}u_{i,j+1} - N_{ij}u_{1,p-2} = 0, \quad i+j = p-2, \quad i,j \ne 0.
\end{equation}
The last relation for $i=2,3, \dots, p-4$ ($i=1$ and $p-3$ give trivial relations) together
with the equality $u_{\frac{p-1}2, \frac{p-1}2} = 0$ (which follows from (2.13)) gives
$p-5$ equations for $p-5$ unknowns $u_{2,p-3}, \dots, u_{p-3,2}$. One easily checks
that
\begin{equation}\notag
u_{ij} = N_{ij}/p \cdot u, \quad i+j = p-1
\end{equation}
for a certain $u \in A$ (actually, $u = -\frac 23 u_{1,p-2}$), provides a unique
solution.

With the aid of (2.11) this equality can be extended to all $i,j, i+j \ge p-1$.
\end{proof}

\begin{lemma}\label{2.6}
Let $d\phi_0 = \xi$, where $\phi_0 \in C_0^2 (W_1(1) \otimes A, W_1(1) \otimes A)$ and 
$\xi \in C_0^3 (W_1(1) \otimes A, W_1(1) \otimes A)$ such that
$\xi (e_i \otimes a, e_j \otimes b, e_k \otimes c)$ is (possibly) nonzero only when
one of the indices $i, j, k$ is equal to $-1$ and the sum of two others is equal to $p-1$.

Then $\xi = 0$ and $\phi_0$ is a cocycle which is cohomologous to $\Upsilon_F + \Psi_D$ 
for some $F \in \mathcal Z^2(A, A)$ and $D \in Der(A)$.
\end{lemma}

\begin{proof}
Write
\begin{equation}\notag
\phi_0 (e_i \otimes a, e_j \otimes b) = \begin{cases}
e_{i+j} \otimes X_{ij}(a,b), & -2 < i+j < p-1 \\
0,                           & \text{otherwise}.
\end{cases}
\end{equation}

Define $\omega \in C^1(W_1(1) \otimes A, W_1(1) \otimes A)$ as follows:
\begin{equation}\notag
\omega (e_{-1} \otimes a) = 0, \quad 
\omega (e_i \otimes a) = e_i \otimes \sum_{j=0}^i X_{-1,j}(1,a), \quad i \ge 0.
\end{equation}
Then 
$d\omega (e_{-1} \otimes 1, e_i \otimes a) = e_{i-1} \otimes X_{-1,i}(1,a) = 
\phi_0 (e_{-1} \otimes 1, e_i \otimes a)$ and replacing
$\phi_0$ by $\phi_0 - d\omega$ (without changing the notation), one can assume that
\begin{equation}\tag{2.14}
X_{-1,j}(1, a) = 0.
\end{equation}

Writing the equation $d\phi_0 = \xi$ for the triple 
$e_{-1} \otimes a, e_{-1} \otimes b, e_i \otimes c$, one obtains
\begin{equation}\tag{2.15}
X_{-1,i-1}(b,ac) - X_{-1,i-1}(a,bc) = aX_{-1,i}(b,c) - bX_{-1,i}(a,c).
\end{equation}
Setting here $b = 1$ and using (2.14), one gets $X_{-1,i-1}(b,c) = X_{-1,i}(b,c)$, 
which implies
\begin{equation}\tag{2.16}
X_{-1,i}(a, b) = X_{-1,0}(a, b).
\end{equation}
Together with (2.15) this gives
\begin{equation}\tag{2.17}
X_{-1,0}(b,ac) - X_{-1,0}(a,bc) - aX_{-1,0}(b,c) + bX_{-1,0}(a,c) = 0.
\end{equation}

Writing the equation $d\phi_0 = \xi$ for the triple 
$e_i \otimes a, e_j \otimes b, e_{-1} \otimes c, i+j \le p-2$, one obtains
\begin{multline}\tag{2.18}
  N_{ij}    X_{-1,i+j}(c, ab) 
-           X_{i-1,j} (ac, b) 
-           X_{i,j-1} (a, bc) 
+           cX_{ij}(a, b)        \\
- N_{i,j-1} aX_{-1,j}(c, b) 
- N_{i-1,j} bX_{-1,i}(c, a) = 0, \quad i+j \le p-2, \quad i,j \ge 0.
\end{multline}
Setting in the last equality $c = 1$, one gets
\begin{equation}\notag
X_{ij}(a,b) = X_{i-1,j}(a,b) + X_{i,j-1}(a,b).
\end{equation}
The last relation together with (2.16) permits to prove, by induction on $i+j$,
the following equality:
\begin{equation}\tag{2.19}
X_{ij}(a,b) = \binom{i+j+1}j X_{-1,0}(a,b) - \binom{i+j+1}i X_{-1,0}(b,a)
\end{equation}

Setting in (2.18) $i=j=0$ and using the fact that 
$X_{00}(a, b) = X_{-1,0}(a, b) - X_{-1,0}(b, a)$ (which follows from (2.19)), 
one obtains
\begin{multline}\tag{2.20}
  X_{-1,0}(bc,a) 
- X_{-1,0}(ac,b) 
- bX_{-1,0}(c,a) 
- cX_{-1,0}(b,a)     \\ 
+ cX_{-1,0}(a,b) 
+ aX_{-1,0}(c,b) = 0.
\end{multline}

Set
\begin{align}
F(a,b) &= \frac 12 (X_{-1,0}(a,b) + X_{-1,0}(b,a) - X_{-1,0}(ab,1))
        = X_{-1,0}(b,a) - aX_{-1,0}(b,1)                            \notag \\
D(a)   &= X_{-1,0}(a,1).                                            \notag
\end{align}

Using (2.17) and (2.20) it is easy to see that $F \in \mathcal Z^2(A, A)$ 
and $D \in Der(A)$, and hence (2.19) implies
\begin{equation}\notag
X_{ij}(a,b) = N_{ij}F(a,b) + \binom{i+j+1}j bD(a) - \binom{i+j+1}i aD(b).
\end{equation}

Thus $\phi_0$ is a cocycle, whence $\xi = 0$.
\end{proof}

\begin{lemma}\label{2.7}
Let $d\phi_p = \xi$, where $\phi_p \in C_p^2(W_1(1) \otimes A, W_1(1) \otimes A)$
and $\xi \in C_p^3(W_1(1) \otimes A, W_1(1) \otimes A)$ 
such that the only possibly nonzero values of $\xi$ are given by
\begin{equation}\notag
\xi (e_{-1} \otimes a, e_{-1} \otimes b, e_0 \otimes c) = e_{p-2} \otimes (aG(b,c) - bG(a,c))
\end{equation}
for some $G \in Hom(A \otimes A, A)$.

Then $G$ is a Harrison 2-coboundary and $\phi_p = \Phi_D$ for some $D \in End(A)$. 
If $G = 0$, then $D \in Der(A)$.
\end{lemma}

\begin{proof}
Write
\begin{equation}\notag
\phi_p (e_i \otimes a, e_j \otimes b) = \begin{cases}
e_{p-2} \otimes X(a,b), & i = j = -1 \\
0,                      & \text{otherwise}.
\end{cases}
\end{equation}
Obviously $X$ is skew-symmetric. Writing the equation $d\phi_p = \xi$ for the triples
$e_{-1} \otimes a, e_{-1} \otimes b, e_{-1} \otimes 1$ and 
$e_{-1} \otimes a, e_{-1} \otimes 1, e_0 \otimes b$, one gets respectively:
\begin{equation}\tag{2.21}
X(a,b) = aX(b,1) - bX(a,1)
\end{equation}
and
\begin{equation}\tag{2.22}
- X(ab, 1) + X(b, a) + 2bX(a, 1) = aG(1, b) - G(a, b).
\end{equation}

Setting $D(a) = X(1,a)$, we obtain $\phi_p = \Phi_D$. 
Substitution of (2.21) into (2.22) gives
\begin{equation}\notag
G(a, b) - aG(1, b) = \delta D(a, b).
\end{equation}
Symmetrizing the last equality, one gets
\begin{equation}\notag
G(a, b) = \delta D(a, b) + abG(1, 1) = \delta (D + R_{G(1,1)})(a, b).
\end{equation}
If $G = 0$ then $\delta D = 0$, i.e. $D \in Der(A)$.
\end{proof}

This completes the proof of Proposition 2.1(ii).

Similar but more elementary computations can be utilized to prove

\begin{proposition}\label{2.8}
\hfill
\begin{enumerate}
\item $H^1 (sl(2) \otimes A, sl(2) \otimes A) \simeq Der(A)$
\item $H^2 (sl(2) \otimes A, sl(2) \otimes A) \simeq Har^2(A, A)$.
\end{enumerate}
\end{proposition}

\begin{proof}
We refer to the paper of Cathelineau \cite{C}. Though formally it contains a
slightly different result -- namely, the computation of 
$H_2 (\mathfrak g \otimes A, \mathfrak g \otimes A)$ for classical simple
Lie algebra $\mathfrak g$ over a field of characteristic zero, 
the methods employed there can be easily
adapted to our case. Alternatively, one may go along the lines of our proof for the
case $W_1(1) \otimes A$. All basic cocycles turn out to be of the type (2.3).
\end{proof}

\section{Low-dimensional cohomology of $\mathfrak L(A, D)$}

\begin{theorem}\label{3.1}
\hfill
\begin{enumerate}
\item $H^2 (\mathfrak L(A, D), K) \simeq (A^*)^D$
\item $H^1 (\mathfrak L(A, D), \mathfrak L(A, D)) \simeq Der(A)^D$
\item $H^2 (\mathfrak L(A, D), \mathfrak L(A, D)) \simeq 
      A^D \oplus Der(A)_D \oplus Der(A)^D \oplus Har^2(A, A)^D$.
\end{enumerate}
\end{theorem}

All super- and subscripts here denote the kernel and cokernel respectively 
of the corresponding action of $D$ 
(which is, for (i), given by $Df(a) \mapsto f(D(a))$ for $f\in A^*$,
and for (ii) and (iii) is the standard action on Harrison (=Hochschild) cocycles 
described in \S 1).

Part (i) borrowed from \cite{me-trans}, where it is proved along the lines of the present paper
(though the computations are easier).

We will give also an explicit basis of $H^1 (\mathfrak L(A, D), \mathfrak L(A, D))$
                                   and $H^2 (\mathfrak L(A, D), \mathfrak L(A, D))$.

There are at least three ways to compute the cohomology of deformed algebra
knowing the cohomology of an initial one. The first way is the Coffee-Gerstenhaber
lifting theory (cf. \cite{GS}) which tells how to determine obstructions to lifting of
cocycles on a Lie algebra $L$ to its deformation $\mathfrak L$.

The second way is applicable when $\mathfrak L$ is a filtered deformation of $L$, 
i.e. $\mathfrak L$ is a filtered Lie algebra with descending filtration 
$\{ \mathfrak L_i \}$ and $L = gr \mathfrak L$. One can define a
descending filtration in the Chevalley-Eilenberg complex $C^*(\mathfrak L, \mathfrak L)$:
\begin{equation}\notag
C_i^n (\mathfrak L, \mathfrak L) = 
\set{\phi \in C^n (\mathfrak L, \mathfrak L)}
{\phi (\mathfrak L_{i_1}, \dots, \mathfrak L_{i_n}) \subseteq \mathfrak L_{i_1 + \dots i_n + i}}.
\end{equation}
Then the associated graded complex will be $C^*(L, L)$ with the grading defined by
(1.4), and the general theory about filtered complexes says that there is a spectral
sequence abutting to $H^*(\mathfrak L, \mathfrak L)$ whose $E_1$-term is $H^*(L,L)$.

The third way is applicable in a special situation when $\mathfrak L$ is
a "1-step" deformation
of $L$, i.e. multiplication in $\mathfrak L$ is given by
\begin{equation}\notag
\{x,y\} = [x,y] + \phi (x,y)t
\end{equation}
where $\liebrack$ is a multiplication in $L$ and $\phi \in Z^2(L, L)$. 
Then we have three complexes defined on the underlying module $C^*(L, L)$: 
the first one responsible for the cohomology of $L$ with differential $d$, 
the second one -- with differential $b = [\,\cdot \,, \phi \,]$ (Massey bracket), and
the third one is responsible for the cohomology of $\mathfrak L$ with differential $b + d$. 
Moreover, the Jacobi identity for $\curlybrack$ implies $bd + db = 0$. 
In this situation it is possible to define a double complex on $C^*(L, L)$ 
whose horizontal arrows are $d$ and vertical ones are $b$. 
The total complex $\mathcal T$ of this double complex is closely related to the 
Chevalley--Eilenberg complex $\mathcal C = C^*(\mathfrak L, \mathfrak L)$ 
responsible for the cohomology of $\mathfrak L$. Namely, there is a surjection
\begin{equation}\notag
\mathcal T^n = \bigoplus_{i=1}^n C^n(L, L) \to C^n(\mathfrak L, \mathfrak L)
\end{equation}
defined by the summation of all coordinates, whose kernel $\mathcal K$ is closely related to
the shifted complex $\mathcal T[-1]$. So one can determine the cohomology 
$H^n(\mathfrak L, \mathfrak L)$ from the long exact sequence associated 
with the short exact sequence of complexes 
$0 \to \mathcal K \to \mathcal T \to \mathcal C \to 0$.

However, in our even more specific situation we will use the fourth method 
employing the special $\mathbb Z$-grading. Its advantage is that we will be able not only to
determine $H^1(\mathfrak L(A, D), \mathfrak L(A, D))$ and 
          $H^2(\mathfrak L(A, D), \mathfrak L(A, D))$ as modules,
but also to find explicit expressions for cocycles.

As noted in \S 1, when considering the cohomology both of $W_1(1) \otimes A$
and $\mathfrak L(A, D)$, we may restrict our attention to a subcomplex preserving the 
$\mathbb Z_p$-grading of $W_1(1) \otimes A$:
\begin{equation}\notag
C_{[0]}^n (W_1(1) \otimes A, W_1(1) \otimes A) = 
\bigoplus_{i \in \mathbb Z} C_{ip}^n (W_1(1) \otimes A, W_1(1) \otimes A).
\end{equation}

Let $d$ and $d_D$ be the differentials in the Chevalley-Eilenberg complexes
$C^*(W_1(1) \otimes A, W_1(1) \otimes A)$ and $C^*(\mathfrak L(A, D), \mathfrak L(A, D)$, 
respectively. We obviously have
\begin{equation}\tag{3.1}
d_D = d + [\, \cdot \,, \Phi_D]
\end{equation}
where $\liebrack$ denotes the graded Lie (super)algebra structure (Massey brackets) on
$H^*(W_1(1) \otimes A, W_1(1) \otimes A)$.

Since $\Phi_D \in C_p^2(W_1(1) \otimes A, W_1(1) \otimes A)$, the bracket 
$b = [\, \cdot \,, \Phi_D]$ acts as a differential of bidegree $(1,p)$ on the 
bigraded module
$C_*^* (W_1(1) \otimes A, W_1(1) \otimes A)$ 
(the first grading is the usual cohomology grading, the second one comes from the 
$\mathbb Z$-grading on $W_1(1) \otimes A$). Denoting for convenience the module
$C_{ip}^n (W_1(1) \otimes A, W_1(1) \otimes A)$ as $\widehat C_i^n$, 
we have a double complex

\def\C{\widehat C}
$$\begin{CD}
\dots   @.     \dots   @.     \dots      @.     \dots                   \\
@AAA           @AAA           @AAA              @AAA                    \\
\C_1^0  @>d>>  \C_1^1  @>d>>  \C_1^2     @>d>>  \C_1^3     @>>>  \dots  \\
@.             @AAbA          @AAbA             @AAbA                   \\       
        @.     \C_0^0  @>d>>  \C_0^1     @>d>>  \C_0^2     @>>>  \dots  \\
@.             @.             @AAbA             @AAbA                   \\
        @.             @.     \C_{-1}^0  @>d>>  \C_{-1}^1  @>>>  \dots  \\
@.             @.             @.                @AAA                    \\
        @.             @.                @.     \dots
\end{CD}$$
\smallskip

In view of (3.1), the total complex of this double complex is exactly the 
Chevalley-Eilenberg complex computing the cohomology 
$H^*(\mathfrak L(A, D), \mathfrak L(A, D))$. Therefore the first spectral sequence
$\{ E_r^{st} \}$ associated with it has the $E_1$-term
\begin{equation}\notag
E_1^{st} \simeq H_{ps}^{s+t}(W_1(1) \otimes A, W_1(1) \otimes A).
\end{equation}

A necessary condition for $\widehat C_i^n \ne 0$ is that there exists a solution
to $-1 \le i_1 + \dots + i_n + ip \le p-2$ for $-1 \le i_k \le p-2$.
This implies the inequalities $-n + \frac{2n-1}p \le i \le 1 + \frac{n-2}p$
(so for $n=1$, $i=0$, for $n=2$, $-1 \le i \le 1$, for $n=3$, $-2 \le i \le 1$).
Thus in each degree there is finite number of nonvanishing components 
and the spectral sequence converges to $H^{s+t}(\mathfrak L(A, D), \mathfrak L(A, D))$.

The only possibly nonzero terms responsible for the cohomology of low degree are:
\begin{equation}\notag
E_r^{01}, E_r^{-1,3}, E_r^{02}, E_r^{11}, E_r^{-2,5}, E_r^{-1,4}, E_r^{03}, E_r^{12}.
\end{equation}

Hence the only possibly nonzero differentials affecting the values of 
$H^1(\mathfrak L(A, D), \mathfrak L(A, D))$ and $H^2(\mathfrak L(A, D), \mathfrak L(A, D))$ are:
\begin{align}
d_1^{01}   &: E_1^{01}   \to E_1^{11}  \notag \\
d_1^{-1,3} &: E_1^{-1,3} \to E_1^{03}  \notag \\
d_1^{02}   &: E_1^{02}   \to E_1^{12}  \notag \\
d_2^{-1,3} &: E_2^{-1,3} \to E_2^{12}.  \notag
\end{align}

Consequently,
\begin{align}
& E_\infty^{01}     = E_1^{01}   = Ker\,d_1^{01}  \notag \\
& E_{\infty}^{-1,3} = E_3^{-1,3} = Ker\,d_2^{-1,3}; \quad E_2^{-1,3} = Ker\,d_1^{-1,3} \notag \\
& E_{\infty}^{02}   = E_2^{02}   = Ker\,d_1^{02}         \tag{3.2}  \\  
& E_{\infty}^{11}   = E_2^{11}   = E_1^{11}/Im\,d_1^{01} \notag
\end{align}
and
$H^1(\mathfrak L(A, D), \mathfrak L(A, D)) \simeq E_\infty^{01}$,
$H^2(\mathfrak L(A, D), \mathfrak L(A, D)) \simeq E_{\infty}^{-1,3} \oplus E_{\infty}^{02}
                                           \oplus E_{\infty}^{11}$.
Proposition 2.1 (strictly speaking, the explicit basic cocycles provided in its proof) yields
\begin{align}
E_1^{01}   &\simeq H_0^1 (W_1(1) \otimes A, W_1(1) \otimes A) \simeq Der(A)  \notag \\
E_1^{-1,3} &\simeq H_{-p}^2 (W_1(1) \otimes A, W_1(1) \otimes A) 
            \simeq H^2(W_1(1), W_1(1)) \otimes A \simeq A                    \notag \\
E_1^{02}   &\simeq H_0^2( W_1(1) \otimes A, W_1(1) \otimes A) 
            \simeq Der(A) \oplus Har^2(A, A)                                 \notag \\
E_1^{11}   &\simeq H_p^2 (W_1(1) \otimes A, W_1(1) \otimes A) \simeq Der(A). \notag
\end{align}

In the next lemmas we will determine all necessary kernels and images in (3.2).

\begin{lemma}\label{3.2}
$E_{\infty}^{-1,3} \simeq A^D$.
\end{lemma}

\begin{proof}
In order to determine $Ker\,d_1^{-1,3}$, one needs to consider the equation
\begin{equation}\tag{3.3}
[(1 \otimes R_u) \circ \Theta_{\phi}, \Phi_D] = d\Lambda_u
\end{equation}
for some $\Lambda_u \in C_0^2 (W_1(1) \otimes A, W_1(1) \otimes A)$.
The terms of $E_2^{-1,3}$ will be of the form 
$(1 \otimes R_u) \circ \Theta_{\phi} - \Lambda_u$ for appropriate solutions of (3.3).

Direct computations show
\begin{equation}\tag{3.4}
[(1 \otimes R_u) \circ \Theta_{\phi}, \Phi_D] = 
(1 \otimes R_u) \circ [\Theta_{\phi}, \Phi_D] + (1 \otimes R_{D(u)}) \circ \Gamma
\end{equation}
where $\Gamma \in C^3(W_1(1) \otimes A, W_1(1) \otimes A)$ defined as 
(assuming $i \le j \le k$)
\begin{equation}\notag
\Gamma (e_i \otimes a, e_j \otimes b, e_k \otimes c) = \begin{cases}
e_{p-2} \otimes N_{jk}/p \cdot abc, & i=-1, j+k=p-1    \\
0,                                  & \text{otherwise}
\end{cases}
\end{equation}
and
\begin{multline}\tag{3.5}
[\Theta_{\phi}, \Phi_D](e_i \otimes a, e_j \otimes b, e_k \otimes c) \\ 
= \begin{cases}
e_{p-2} \otimes N_{jk}/p \, (bcD(a) - aD(bc)), & i = -1, j + k = p-1  \\
0,                                          & \text{otherwise}.
\end{cases}
\end{multline}
(Notice that $N_{jk}/p = (-1)^k \frac{2k+1}{k(k+1)}$ if $j + k = p-1$).

Define
\begin{equation}\notag
\Theta^{\prime} (e_i \otimes a, e_j \otimes b) = \begin{cases}
e_{i+j} \otimes (\lambda_{ij} aD(b) - \lambda_{ji} bD(a)), & -2 < i + j < p-1  \\
0,                                                         & \text{otherwise}
\end{cases}
\end{equation}
for some $\lambda_{ij} \in K$.

Let us verify first that there are $\lambda_{ij}$ such that 
$\Theta_{\phi} + \Theta^{\prime} \in E_2^{-1,3}$. Writing the
equation (3.3) under conditions $u = 1$ and $\Lambda_u = -\Theta^{\prime}$ 
for the triple $e_{-1} \otimes a, e_j \otimes b, e_k \otimes c, j + k = p-1$
and for all remaining cases, we obtain respectively:
\begin{align}
\lambda_{j-1,k} + \lambda_{j,k-1} &= (-1)^k \frac{2k+1}{k(k+1)}   \tag{3.6}  \\
\lambda_{j,k-1} - \lambda_{k,j-1} &= 
2(-1)^k \lambda_{j,-1} + 2(-1)^{k+1} \lambda_{k,-1} + (-1)^{k+1} \frac{2k+1}{k(k+1)} \tag{3.7}
\end{align}
where $j + k = p-1$, and
\begin{multline}\tag{3.8}
  N_{ij}    \lambda_{i+j,k} 
- N_{jk}    \lambda_{i,j+k} 
+ N_{ik}    \lambda_{j,i+k} 
+ N_{j+k,i} \lambda_{jk} 
- N_{i+k,j} \lambda_{ik} \\ 
= 0, \quad i,j,k \ge 0, \quad i+j+k<p-1.
\end{multline}
(the left-hand side in the latter is a basic expression for the coefficient of $abD(c)$ in
$d\Theta^{\prime}$).

\begin{lemma}\label{3.3}
\begin{equation}\notag
\lambda_{ij} = \sum_{k=1}^i \binom{i+j+1-k}{j+1} \frac{k+2}{k(k+1)}, \quad -1 \le i,j \le p-2
\end{equation}
provides solution for (3.6)--(3.8).
\end{lemma}

\proof
Note the following properties of the just defined coefficients $\lambda_{ij}$:
\begin{enumerate}
\item $\lambda_{-1,j} = \lambda_{0j} = 0; \quad \lambda_{1,j} = \frac 32$
\item $\lambda_{i,-1} = \sum_{k=1}^i \frac{k+2}{k(k+1)}$
\item $\lambda_{ij} = \lambda_{i-1,j} + \lambda_{i,j-1}$.
\end{enumerate}

Now, (3.6) may be reformulated as
\begin{equation}\notag
\lambda_{ij} = (-1)^j \frac{2j+1}{j(j+1)}, \quad i+j=p-1, \quad i,j \ge 1
\end{equation}
which can be proved with the help of simple transformations of binomial coefficients
in the spirit of the first few pages of \cite{riordan}.

(3.7) is proved by induction on $j$, using (3.6) in the induction step.

Finally, (3.8) is proved by induction on $i + j + k$. The induction step is:
\begin{multline}\notag
  N_{ij}    \lambda_{i+j,k} 
- N_{jk}    \lambda_{i,j+k} 
+ N_{ik}    \lambda_{j,i+k} 
+ N_{j+k,i} \lambda_{jk} 
- N_{i+k,j} \lambda_{ik}          \\
=                                 
  N_{ij}      \lambda_{i+j-1,k} 
+ N_{ij}      \lambda_{i+j,k-1} 
- N_{jk}      \lambda_{i-1,j+k} 
- N_{jk}      \lambda_{i,j+k-1} 
+ N_{ik}      \lambda_{j-1,i+k}   \\
+ N_{ik}      \lambda_{j,i+k-1} 
+ N_{j+k,i-1} \lambda_{jk}
+ N_{j+k-1,i} \lambda_{jk} 
- N_{i+k,j-1} \lambda_{ik} 
- N_{i+k-1,j} \lambda_{ik}        \\
= 
(
  N_{i-1,j}   \lambda_{i+j-1,k}
- N_{jk}      \lambda_{i-1,j+k}
+ N_{i-1,k}   \lambda_{j,i+k-1} 
+ N_{j+k,i-1} \lambda_{jk} 
- N_{i+k-1,j} \lambda_{i-1,k}
)                                 \\ 
+ (
  N_{i,j-1}   \lambda_{i+j-1,k}
- N_{j-1,k}   \lambda_{i,j+k-1}
+ N_{ik}      \lambda_{j-1,i+k} 
+ N_{j+k-1,i} \lambda_{j-1,k} 
- N_{i+k,j-1} \lambda_{ik}
)                                 \\ 
+ (
  N_{ij}      \lambda_{i+j,k-1} 
- N_{j,k-1}   \lambda_{i,j+k-1}
+ N_{i,k-1}   \lambda_{j,i+k-1}
+ N_{j+k-1,i} \lambda_{j,k-1}
- N_{i+k-1,j} \lambda_{i,k-1}
)                                 \\ 
= 0
\end{multline}
where the first equality follows from recurrent relations for $\lambda_{ij}$, 
the second one from those for $N_{ij}$, and the third one from the induction assumption 
for triples $(i-1,j,k)$, $(i,j-1,k)$ and $(i,j,k-1)$.
\end{proof}

\noindent \textit{Continuation of the proof of Lemma 3.2}. 
Now consider a general solution of (3.3).
Taking into account (3.4), the partial solution 
$[\Theta_{\phi}, \Phi_D] = -d\Theta^{\prime}$, and the commutativity of operators 
$d$ and $R_u$, (3.3) can be rewritten as
\begin{equation}\notag
d(\Lambda_u + (1 \otimes R_u) \circ \Theta^{\prime}) = (1 \otimes R_{D(u)}) \circ \Gamma.
\end{equation}

By Lemma 2.6, $D(u) = 0$ and $\Lambda_u = - (1 \otimes R_u) \circ \Theta^{\prime}$ 
up to elements from $Z_0^2 (W_1(1) \otimes A, W_1(1) \otimes A)$. 
Hence $E_2^{-1,3}$ consists of elements of the form
\begin{equation}\notag
\overline \Theta_u = (1 \otimes R_u) \circ (\Theta_{\phi} + \Theta^{\prime}), \quad u \in A^D
\end{equation}
and $E_2^{-1,3} \simeq A^D$.

To compute $Ker\,d_2^{-1,3}$, take a look at $Im\,d_2^{-1,3}$, i.e. on elements of the form
$[(1 \otimes R_u) \circ (\Theta_{\phi} + \Theta^{\prime}), \Phi_D], u \in A^D$. 
The latter expression is equal to $(1 \otimes R_u) \circ [\Theta^{\prime}, \Phi_D]$
up to elements from $B^3(W_1(1) \otimes A, W_1(1) \otimes A)$. Direct computations show
\begin{multline}\notag
[\Theta^{\prime}, \Phi_D] (e_i \otimes a, e_j \otimes b, e_k \otimes c) \\
= \begin{cases}
e_{p-2} \otimes \lambda_{p-2,0}(aD(b) - bD(a))D(c), & i = j = -1, k = 0  \\
0,                                                  & \text{otherwise}.
\end{cases}
\end{multline}

But 
\begin{equation}\notag
\lambda_{p-2,0} = - \sum_{k=1}^{p-2} (1 + \frac 2k) = 
-(p-2) - 2\sum_{k=1}^{p-1} k + \frac 2{p-1} = 0.
\end{equation}

Consequently, $d_2^{-1,3}$ is zero and 
$E_{\infty}^{-1,3} = E_3^{-1,3} = E_2^{-1,3} \simeq A^D$. \qed

\begin{lemma}\label{3.4}
$E_{\infty}^{02} \simeq Der(A)^D \oplus Har^2(A,A)^D$.
\end{lemma}

\begin{proof}
To determine $Ker\,d_p^{02}$, one needs to solve two equations
\begin{align}
[\Psi_E, \Phi_D]     &= d\Lambda_E     \tag{3.9}  \\
[\Upsilon_F, \Phi_D] &= d\Lambda_F     \tag{3.10}
\end{align}
for some $\Lambda_E, \Lambda_F \in C_p^2 (W_1(1) \otimes A, W_1(1) \otimes A)$. Let
\begin{equation}\notag
\Psi_E^{\prime} (e_i \otimes a, e_j \otimes b) = \begin{cases}
e_{p-2} \otimes (D(a)E(b) - E(a)D(b)), & i=j=-1           \\
0,                                     & \text{otherwise}.
\end{cases}
\end{equation}

By means of direct computations one gets
\begin{equation}\notag
[\Psi_E, \Phi_D] = d\Psi_E^{\prime} + \Gamma_E
\end{equation}
where the only possibly nonzero values of $\Gamma_E$ are given by
\begin{equation}\notag
\Gamma_E (e_{-1} \otimes a, e_{-1} \otimes b, e_{-1} \otimes c) = 
e_{p-2} \otimes (a[E,D](b) - b[E,D](a))c.
\end{equation}

But (3.9) implies that $\Gamma_E$ is a coboundary $d(\Lambda_E - \Psi_E^{\prime})$. 
One easily checks that each coboundary vanishes on the triple 
$e_{-1} \otimes a, e_{-1} \otimes 1, e_0 \otimes 1$, which implies
$[E,D] = 0$ and $\Gamma_E = 0$. Consequently, $\Lambda_E = \Psi_E^{\prime}$ 
up to elements from $Z_p^2 (W_1(1) \otimes A, W_1(1) \otimes A)$ and the set of elements
$\set{\Psi_E - \Psi_E^{\prime}}{E \in Der(A)^D}$ embeds into $E_2^{02}$.

To solve equation (3.10), define
\begin{equation}\notag
\Upsilon_F^{\prime} (e_i \otimes a, e_j \otimes b) = \begin{cases}
e_{p-2} \otimes (F(D(a),b) - F(a,D(b))), & i=j=-1  \\
0,                                       & \text{otherwise}.
\end{cases}
\end{equation}

By means of direct computations one gets
\begin{equation}\tag{3.11}
[\Upsilon_F, \Phi_D] = d\Upsilon_F^{\prime} + \Gamma_F
\end{equation}
where the only possibly nonzero values of $\Gamma_F$ are given by
\begin{equation}\notag
\Gamma_F (e_{-1} \otimes a, e_{-1} \otimes b, e_0 \otimes c) = 
e_{p-2} \otimes (aD \star F(b,c) - bD \star F(a, c)).
\end{equation}

By (3.11), $\Gamma_F = d(\Lambda_F - \Upsilon_F^{\prime})$. 
According to Lemma 2.7, $\Lambda_F - \Upsilon_F^{\prime} = \Phi_H$ for some
$H \in End(A)$ and moreover, $D \star F = \delta H$ 
(we may suppose that $D \star F(1,1) = 0$).

Conversely, if $D \star F = \delta H$, then $\Gamma_F = -d\Phi_H$,
which in view of (3.11) leads to the solution 
$\Lambda_F = \Upsilon_F^{\prime} + \Phi_H$ of (3.10).

So, $E_{\infty}^{02} = E_2^{02}$ is the direct sum of two subspaces consisting of elements
of the form $\Psi_E - \Psi_E^{\prime}$ and $\Upsilon_F - \Upsilon_F^{\prime} - \Phi_H$
for appropriate $E$, $F$ and $H$, and isomorphic to $Der(A)^D$ and $Har^2(A,A)^D$ 
respectively. 
\end{proof}

\begin{lemma}\label{3.5}
\hfill
\begin{enumerate}
\item $E_\infty^{01}   \simeq Der(A)^D$
\item $E_{\infty}^{11} \simeq Der(A)_D$.
\end{enumerate}
\end{lemma}

\begin{proof}
$d_1^{01}$ acts on the space $E_1^{01} \simeq Der(A)$ as
$1\otimes E \mapsto [1\otimes E, \Phi_D] = \Phi_{[D,E]}$, $E\in Der(A)$.
Hence $Ker\,d_1^{01} \simeq Der(A)^D$, proving (i).

$Im\,d_1^{01} \simeq [D, Der(A)]$, 
$E_{\infty}^{11} = E_2^{11} \simeq Der(A)_D$ and $E_{\infty}^{11}$ consists
of elements $\Phi_E$ for $E \in Der(A)$ which are independent modulo $[D, Der(A)]$,
proving (ii).
\end{proof}

Putting all these calculations together, we get statements (ii) and (iii) of Theorem 3.1.
(Lemma 3.5(i) is used to get a formula for cohomology of degree 1,
while all the rest is used to get a formula for cohomology of degree 2).

For convenience we summarize here the cocycles whose cohomology
classes constitute a basis of $H^1(\mathfrak L(A, D), \mathfrak L(A, D))$
                          and $H^2(\mathfrak L(A, D), \mathfrak L(A, D))$.

Basic cocycles of degree 1 are just mappings of the form $1\otimes E$, $E\in Der(A)^D$.

All cocycles of degree 2 constructed here have their counterparts in 
$Z^2(W_1(1) \otimes A, W_1(1) \otimes A)$ (in fact, they are liftings, in
the Gerstenhaber's terminology \cite{GS}, of 2-cocycles on $W_1(1) \otimes A$). 
Each class of cocycles denoted by overlined capital Greek letter is lifted from the 
corresponding class of \S 2 denoted by the same letter.

So, let $\overline \Theta_u$, $\overline \Upsilon_{F,H}$, $\overline \Psi_E$ and 
$\overline \Phi_E$ be 2-cochains on $\mathfrak L(A, D)$ defined by the following
formulas, where the top line comes from the appropriate cocycle of \S 2
(the ``regular'' components), and the second line represent a new component coming from
the deformation:
\begin{equation}\notag
\overline \Theta_u (e_i \otimes a, e_j \otimes b) = \begin{cases}
e_{i+j-p} \otimes N_{ij}/p \> abu,                            & i+j \ge p-1      \\
e_{i+j}   \otimes (\lambda_{ij} aD(b) - \lambda_{ji} bD(a))u, & -2 < i+j < p- 1  \\
0,                                                            & \text{otherwise}
\end{cases}
\end{equation}
where $u \in A^D$ and the coefficients $\lambda_{ij}$ defined as in Lemma 3.3
(the regular component is (2.2)),
\begin{multline}\notag
\overline \Upsilon_{F,H} (e_i \otimes a, e_j \otimes b)  \\
= \begin{cases}
e_{i+j} \otimes N_{ij} F(a,b),                           & -2<i+j<p-1  \\
e_{p-2} \otimes (bH(a) - aH(b) - F(D(a),b) + F(a,D(b))), & i = j = -1  \\
0,                                                       & i+j \ge p-1
\end{cases}
\end{multline}
where $F \in \mathcal Z^2(A, A)^D$ and $H \in End(A)$ such that $D \star F = \delta H$
(the regular component is (2.3)),
\begin{equation}\notag
\overline \Psi_E (e_i \otimes a, e_j \otimes b) = \begin{cases}
e_{p-2} \otimes (E(a)D(b) - E(b)D(a)),                         & i = j = -1        \\
e_{i+j} \otimes (\binom{i+j+1}j bE(a) - \binom{i+j+1}i aE(b)), & -2 < i + j < p-1  \\
0,                                                             & i+j \ge p-1
\end{cases}
\end{equation}
where $E \in Der(A)^D$ (the regular component is (2.5)), 
and, finally, $\overline \Phi_E = \Phi_E$ (the regular component is (1.2); 
there is no deformation component).

Lemma 2.4 (stating the independence of initial cocycles on $W_1(1) \otimes A$) together
with the spectral sequence construction assure the independence of the corresponding
cocycles on $\mathfrak L(A, D)$. More precisely, the following is true:

\begin{proposition}\label{3.6}
Let $\{ u_i \}$ be linearly independent elements of $A$, $\{ F_i \}$ be 
cohomologically independent cocycles in $\mathcal Z^2(A,A)$, $\{ H_i \}$ be 
elements in $Hom(A,A)$ linearly independent modulo $Der(A)$ and such that 
$D \star F_i = \delta H_i$, $\{ E_i \}$ be linearly independent elements in $Z_D(Der(A))$,
and $\{ E_i^{\prime} \}$ be derivations of $A$ linearly independent modulo $[D, Der(A)]$.

Then the cocycles $\overline \Theta_{u_i}$, $\overline \Upsilon_{F_i,H_i}$, 
$\overline \Psi_{E_i}$, $\overline \Phi_{E_i^{\prime}}$ are cohomologically independent. 
\end{proposition}

\section{Low-dimensional cohomology of $W_1(n) \otimes A$}

Now our objective is to transform the results obtained so far for $\mathfrak L(A, D)$ into
those for $W_1(n) \otimes A$.

For this, take $A = O_1(n-1) \otimes B$ and $D = \partial \otimes 1$. 
By Proposition 1.2, $\mathfrak L(A, D)$ in this case isomorphic to $W_1(n) \otimes B$, 
and Theorem 3.1 entails
\begin{multline}\tag{4.1}
H^2(W_1(n) \otimes B, W_1(n) \otimes B)                \\ 
\simeq (O_1(n-1) \otimes B)^{\partial \otimes 1} 
\oplus Der(O_1(n-1) \otimes B)_{\partial \otimes 1} 
\oplus Der(O_1(n-1) \otimes B)^{\partial \otimes 1}    \\
\oplus Har^2(O_1(n-1) \otimes B, O_1(n-1) \otimes B)^{\partial \otimes 1}.
\end{multline}

The next lemmas collect all necessary information for evaluation of four summands 
appearing on the right side of this isomorphism. 
(Just for notational convenience, we put $m = n-1$).

\begin{lemma}\label{4.1}
\hfill
\begin{enumerate}
\item $(O_1(m) \otimes B)^{\partial \otimes 1} = 1 \otimes B$
\item $Der(O_1(m) \otimes B)_{\partial \otimes 1} \simeq 
      \langle x^{p^m-1}\partial^{p^k} \,|\, 0 \le k \le m-1 \rangle \otimes B \>\oplus\>
      x^{p^m-1} \otimes Der(B)$
\item $Der(O_1(m) \otimes B)^{\partial \otimes 1} \simeq
      \langle \partial^{p^k} \,|\, 0 \le k \le m-1 \rangle \otimes B \>\oplus\>
      1 \otimes Der(B)$.
\end{enumerate}
\end{lemma}

\begin{proof}
(i) Obvious, as $Ker_{O_1(m)}\partial = K1$.

(ii) Since
\begin{equation}\notag
Der(O_1(m) \otimes B) \simeq Der(O_1(m)) \otimes B + O_1(m) \otimes Der(B)
\end{equation}
and $Der(O_1(m))$ is a free $O_1(m)$-module with basis 
$\set{\partial^{p^k}}{0 \le k \le m-1}$,
\begin{multline}\notag
[\partial \otimes 1, Der(O_1(m) \otimes B)] \simeq 
[\partial, Der(O_1(m))] \otimes B \>\oplus\> \partial(O_1(m)) \otimes Der(B)   \\ 
= \langle x^i \partial^{p^k} \,|\, 0 \le i < p^m-1, 0 \le k \le m-1 \rangle \otimes B 
\>\oplus\>
\langle x^i \,|\, 0 \le i < p^m-1 \rangle \otimes Der(B).
\end{multline}

As $\langle x^{p^m-1}\partial^{p^k} \,|\, 0 \le k < m-1 \rangle$ is a complement in 
$Der(O_1(m))$ to the tensor factor in the first summand, and $\langle x^{p^m-1} \rangle$
is a complement in $O_1(m)$ to those in the second summand, we get the isomorphism 
desired.

(iii) Analogous to (ii).
\end{proof}

Further, according to \cite{harrison}, Theorem 5,
\begin{multline}\tag{4.2}
Har^2(O_1(m) \otimes B, O_1(m) \otimes B)^{\partial \otimes 1}   \\
\simeq Har^2(O_1(m),O_1(m))^{\partial} \otimes B \>\oplus\>
 O_1(m)^{\partial} \otimes Har^2(B,B)
\end{multline}
(as $O_1(m)^{\partial} \simeq K$, the second summand is actually just $Har^2(B,B)$).

So we need to compute the second Harrison cohomology of the divided powers algebra 
$O_1(m)$. First we determine its Hochschild cohomology.
It is more convenient to work with reduced polynomial ring $O_m$. 

Note that $O_m$ is a factor-algebra of a polynomial algebra as well as the
group algebra of an elementary abelian group, and for both class of algebras all sort of
cohomological computations have been done. Instead of digging the result
we need out of the literature (which will require some additional computations anyway,
see e.g. \cite{loday}, \S 7.4 and \cite{holm} and references therein), we 
give a direct simple proof suited for our needs.

We use multi-index notations:
$\digamma_m = 
\set{\alpha = (\alpha_1, \dots, \alpha_m) \in \mathbb Z^m}{0 \le \alpha_i < p}$, 
$x^\alpha = x_1^{\alpha_1} \dots x_m^{\alpha_m}$, $\varepsilon_i$ denotes element
in $\digamma_m$ of the form $(0, \dots, 1 , \dots, 0)$ (1 in the $i$th position).

\begin{proposition}\label{4.2}
$H^i(O_m,O_m)$ is a free $O_m$-module of dimension $\binom{i+m-1}i$.
\end{proposition}

\begin{proof}
Consider first the case $m = 1$, i.e., an algebra $O_1 = K[x]/(x^p)$. 
We use a very simple (and nice) free $O_1^e$-resolution of $O_1$ presented in \cite{RS}
(as $O_1$ is commutative, $O_1^e \simeq O_1 \otimes O_1$):
\begin{equation}\notag
\dots           \overset{d_1}\longrightarrow 
O_1 \otimes O_1 \overset{d_0}\longrightarrow 
O_1 \otimes O_1 \overset{m}\longrightarrow
O_1             \longrightarrow
0
\end{equation}
where $m$ is the multiplication in $O_1$ and
\begin{equation}\notag
d_i (a \otimes b) = \begin{cases}
a \otimes xb - ax \otimes b,              & i \text{ even}  \\
\sum_{k=0}^{p-1} ax^k \otimes x^{p-1-k}b, & i \text{ odd}.
\end{cases}
\end{equation}

Applying the functor $Hom_{O_1^e}(-, O_1)$, we get a (deleted) complex whose all but first
differentials are zero:
\begin{equation}\notag
0   \longrightarrow 
O_1 \overset{id}\longrightarrow 
O_1 \overset{0}\longrightarrow
O_1 \overset{0}\longrightarrow
\dots
\end{equation}

Therefore, for each $i$, $H^i(O_1,O_1) \simeq O_1$.

Now the general case is proved via induction on $m$ by applying the K\"unneth
formula to the decomposition $O_m \simeq O_{m-1} \otimes O_1$.
\end{proof}

\begin{lemma}\label{4.3}
\hfill
\begin{enumerate}
\item $Har^2(O_m,O_m)$ is a free $O_m$-module of dimension $m$. The basic cocycles
(over $O_m$) can be chosen as
\begin{equation}\notag
F_i (x^{\alpha}, x^{\beta}) = \begin{cases}
x^{\alpha + \beta - p\varepsilon_i}, & \alpha_i + \beta_i \ge p  \\
0,                                   & \alpha_i + \beta_i < p
\end{cases}
\end{equation}
for $1 \le i \le m$.
\item $\dim Har^2(O_1(m),O_1(m))^{\partial} = m$. The basic $\partial$-invariant cocycles
can be chosen as
\begin{equation}\tag{4.3}
F_i (x^{\alpha}, x^{\beta}) = \begin{cases}
\binom{\alpha + \beta}{\beta}/p \cdot x^{\alpha + \beta - p^i}, & \alpha_i + \beta_i \ge p \\
0,                                                              & \alpha_i + \beta_i < p
\end{cases}
\end{equation}
where $\alpha = \sum_{i \ge 1} \alpha_i p^{i-1}, \beta = \sum _{i \ge 1} \beta_i p^{i-1}$
are $p$-adic decompositions.
\end{enumerate}
\end{lemma}

\begin{proof}
(i) By Proposition 4.2, $H^2(O_m,O_m)$ is a free $O_m$-module of dimension $\frac{m(m+1)}2$.
We assert that the two classes of cocycles, $F_i, 1 \le i \le m$ and 
$\partial / \partial x_i \cup \partial / \partial x_j, 1 \le i < j \le m$, 
form a basis of this module. Indeed, the cocycle condition is verified
immediately. As we have $m + \frac{m(m-1)}2 = \frac{m(m+1)}2$ cocycles, it remains to check their
independence. Since $F_i$ are symmetric and 
$\partial / \partial x_i \cup \partial / \partial x_j$ are skew, one suffices to
do this only for $F_i$ (remember that 2-coboundaries are symmetric). Suppose
\begin{equation}\tag{4.4}
\sum_{i=1}^m u_i F_i = \delta G
\end{equation}
for some $G \in Hom(O_m,O_m)$ and $u_i \in O_m$.

Then $\delta G(x^{\alpha}, x^{\beta}) = 0$ if $\alpha_i + \beta_i < p$ for each $i$. 
This implies that $G$ acts as derivation on products $x^\alpha x^\beta$ if 
$\alpha_i + \beta_i < p$ for all $i$, hence
\begin{equation}\notag
G(x^{\alpha}) = \sum_{i=1}^m \alpha_i x^{\alpha - \varepsilon_i} G(x^{\varepsilon_i})
\end{equation}
what in its turn entails that $G$ is a derivation, and thus
$\delta G(x^\alpha, x^\beta) = 0$ for all $\alpha, \beta$. 
Then evaluating the left side of (4.4) for all pairs $(x^\alpha, x^\beta)$ such that
$\alpha_j + \beta_j = \delta_{ji}p$ for each $j$ and a fixed $i$, we get $u_i = 0$.
This shows that cocycles $F_i$ are independent.

Now picking from the basic cocycles of $H^2(O_m,O_m)$ those which are symmetric,
we obtain a basis $\set{F_i}{1 \le i \le m}$ of $Har^2(O_m,O_m)$ (as a module over $O_m$).
The freeness of $Har^2(O_m,O_m)$ follows either from the previous reasonings or from the
fact that the Harrison cohomology is a direct summand of the Hochschild one (cf. \cite{GS}).

(ii) Using the isomorphism (1.1), the cocycles of part (i) may be rewritten as
(4.3). Direct easy check shows that the cocycles $F_i$ are $\partial$-invariant 
(in fact, $\partial \star F_i = 0$). The identity 
$\partial \star (uF) = u \star \partial F - (\partial u)F$ shows that the equality
$\partial \star (u_1F_1 + \dots + u_kF_k) = 0$ implies
\begin{equation}\notag
(\partial u_1)F_1 + \dots + (\partial u_k)F_k = 0
\end{equation}
which due to the freeness of $Har^2(O_1(m), O_1(m))$ over $O_1(m)$ entails that all 
$u_i \in K1$, and the assertion desired follows.
\end{proof}

Now, collecting (4.1), (4.2) and Lemmas 4.1 and 4.3(ii), we get an isomorphism
\begin{equation}\tag{4.5}
H^2(W_1(n) \otimes B, W_1(n) \otimes B) \simeq 
\mathcal H \otimes B \>\oplus\> Der(B) \>\oplus\> Der(B) \>\oplus\> Har^2(B, B)
\end{equation}
where $\mathcal H$ is a vector space with basis 
$\set{1,x^{p^{n-1}-1} \partial^{p^k} , \partial^{p^k}, F_{k+1}}{0 \le k \le n-2}$.

To obtain an explicit basis of this cohomology group, let us regroup the basis of
$H^2(\mathfrak L(A, D), \mathfrak L(A, D))$, exhibited in \S 3, 
according to the direct summands in (4.5) as follows.

The classes of $3n-2$ cocycles
\begin{equation}\notag
\overline\Theta_{1\otimes u}, 
\overline\Upsilon_{F_{i+1}\otimes R_u,0},
\overline\Psi_{\partial^{p^i}\otimes R_u},
\overline\Phi_{x^{p^{n-1}-1}\partial^{p^i}\otimes R_u}, 
\quad 0 \le i \le n-2, \quad u\in B 
\end{equation}
form a module denoted in (4.5) as $\mathcal H \otimes B$. 
It is easy to see that all these cocycles are of the form 
$(1 \otimes R_u) \circ \Theta_{\phi}$ for appropriate $\phi \in Z^2(W_1(n), W_1(n))$. 
As by Proposition 3.6 all these cocycles are independent, the corresponding $3n-2$ cocycles
on $W_1(n)$ are also independent. But according to \cite{DK}, $\dim H^2(W_1(n),W_1(n)) = 3n-2$,
whence $\mathcal H \simeq H^2(W_1(n), W_1(n))$. It should be noted that the 2-cocycles on $W_1(n)$
derived here do not wholly coincide with basic cocycles presented in \cite{DK}.

The classes of cocycles $\overline\Psi_{1\otimes D}$ and 
$\overline\Phi_{x^{p^{n-1}-1}\otimes D}, D \in Der(B)$, denoted from now for the
convenience as $\psi_D$ and $\phi_D$ respectively, form two modules isomorphic to $Der(B)$.
They are just obvious generalizations of cocycles $\Psi_D$ and $\Phi_D$ to arbitrary $n$:
\begin{multline}\tag{4.6}
\psi_D (e_i \otimes a, e_j \otimes b) = e_{i+j} \otimes 
(\binom{i+j+1}j bD(a) - \binom{i+j+1}i aD(b)),            \\
-1 \le i,j \le p^n - 2
\end{multline}
\begin{equation}\tag{4.7}
\phi_D (e_i \otimes a, e_j \otimes b) = \begin{cases}
e_{p^n-2} \otimes (aD(b) - bD(a)), & i = j = -1  \\
0,                                 & \text{otherwise}.
\end{cases}
\end{equation}

And finally, the classes of cocycles $\overline\Upsilon_{1\otimes F,0}$ where 
$F \in \mathcal Z^2(B,B)$, generate a module isomorphic to $Har^2(B,B)$. 
These cocycles are of the form $\Upsilon_F$ (cf. Proposition 2.2).

Thus we get a generalization of Proposition 2.1:

\begin{theorem}\label{4.4}
For an arbitrary associative commutative unital algebra $B$,
\begin{multline}\notag
H^2(W_1(n)\otimes B, W_1(n)\otimes B)   \\ 
\simeq
H^2(W_1(n), W_1(n)) \otimes B \>\oplus\> Der(B) \>\oplus\> Der(B) \>\oplus\> Har^2(B, B).
\end{multline}
The basic cocycles can be chosen among $(1\otimes R_u)\circ \Theta_\phi$ for 
$\phi\in Z^2(W_1(n),W_1(n))$, $\psi_D$, $\phi_D$ for $D\in Der(B)$, and
$\Upsilon_F for F\in \mathcal Z^2(B,B)$, given by formulas (2.2), (4.6), (4.7) and (2.3) 
respectively.
\end{theorem}

We conclude this section with formulation of all necessary results needed for
our further purposes, which are obtained in a similar (and much simpler) way as
Theorem 4.4 and/or can be found elsewhere (cf. \cite{block-diff}, \cite{C}, \cite{me-trans}):
\begin{align}
H^1 (W_1(n) \otimes B, W_1(n) \otimes B) &\simeq 
H^1(W_1(n),W_1(n)) \otimes B \>\oplus\> Der(B)               \tag{4.8}  \\
H^2 (W_1(n) \otimes B, K) &\simeq H^2(W_1(n),K) \otimes B^*  \tag{4.9}  \\
H^2 (sl(2) \otimes A, K)  &\simeq HC^1(A).                   \tag{4.10}
\end{align}

\section{Filtered deformations of $W_1(n) \otimes A + 1 \otimes \mathfrak D$}

As explained in \S 1, we are interested in filtered Lie algebras whose associated
graded algebra is $S\otimes O_m + 1\otimes \mathfrak D$ for $S = W_1(n)$ or $sl(2)$,
where $\mathfrak D$ is a subalgebra of $Der(O_m)$. First, basing on Theorem 4.4, 
we shall compute the second cohomology group of such algebras.

\begin{lemma}\label{5.1}
Let $L$ be a Lie algebra which can be written as the semidirect product
$L = I \oplus Q$, where $I$ is a centerless perfect ideal of $L$, $Q$ is a subalgebra,
and $Q \cap ad(I) = 0$ (in the last equality, $Q$ and $ad(I)$ considered as 
subspaces of $End(I)$). Then the terms relevant to the cohomology group $H^2(L,L)$
in the Hochschild-Serre spectral sequence of L with respect to I
with the general $E_2$-term $E_2^{pq} = H^p(Q, H^q(I,L))$, are the following:
\begin{align}
E_{\infty}^{20} &= 0                              \notag \\
E_{\infty}^{11} &= E_2^{11} = H^1(Q, H^1(I,I)/Q)  \notag \\
E_2^{02}        &= H^2(I,I)^Q \oplus (Ker\,F)^Q   \notag \\ 
E_{\infty}^{02} &= E_3^{02} = Ker\,d_2^{02}       \notag
\end{align}
where $F: H^2(I) \otimes Q \to H^3(I,I)$ is induced by the mapping
\begin{align}
C^2(I,Q) &\to C^3(I,I)                    \notag \\
\phi     &\mapsto (x \wedge y \wedge z \mapsto [x, \phi(y,z)] + \curvearrowright). \notag
\end{align}
\end{lemma}

\begin{proof}
One has $E_2^{p0} = H^p(Q,L^I)$. The condition $Q \cap ad(I) = 0$ entails 
$L^I = Z(I) = 0$, so $E_2^{p0} = 0$. 
Thus $E_{\infty}^{20} = 0$ and $E_\infty^{02} = E_3^{02}$
follow from standard considerations.
As $d_2$ maps $E_2^{11}$ to $E_2^{30} = 0$, $E_\infty^{11} = E_2^{11}$.
We also have
\begin{align}
E_2^{02} &= H^0(Q,H^2(I,L)) = H^2(I,L)^Q  \notag \\
E_2^{11} &= H^1(Q,H^1(I,L))               \notag \\
E_2^{21} &= H^2(Q,H^1(I,L)).              \notag
\end{align}

Consider a piece of the cohomology long exact sequence associated with the
short exact sequence $0 \to I \to L \to Q \to 0$ of $I$-modules 
($Q$ considered as a trivial $I$-module):
\begin{multline}\tag{5.1}
H^0(I,L) \to H^0(I,Q) \to H^1(I,I) \to H^1(I,L) \to H^1(I,Q) \to H^2(I,I)  \\
\to H^2(I,L) \to H^2(I,Q) \overset{F}\to H^3(I,I)
\end{multline}
($F$ is connecting homomorphism).

We obviously have: $H^0(I,L) = L^I = 0$, $H^0(I,Q) = Q^I = Q$, 
$H^1(I,Q) = H^1(I) \otimes Q = 0$, $H^2(I,Q) = H^2(I) \otimes Q$. Hence
\begin{align}
H^1(I,L) &\simeq H^1(I,I)/Q              \notag \\
H^2(I,L) &\simeq H^2(I,I) \oplus Ker\,F  \notag
\end{align}
(note that since $Q \cap ad(I) = 0$, $Q$ consists of outer derivations of $I$, 
and therefore embeds in $H^1(I,I)$).

As $L = I \oplus Q$ as $Q$-modules and differential commutes with each
$ad\,x, x\in Q$, the $Q$-action commutes with inclusion and projection arrows in (5.1) 
(but not necessarily with connecting homomorphism), and we get
\begin{equation}\notag
H^2(I,L)^Q \simeq H^2(I,I)^Q \oplus (Ker\,F)^Q.
\end{equation}

Putting all this together, we obtain the asserted equalities.
\end{proof}

Passing to our specific case, define a grading on $L = S \otimes B + 1\otimes\mathfrak D$
as in Theorem 1.5, i.e.
\begin{equation}\tag{5.2}
L_i = \begin{cases}
e_0 \otimes B + 1\otimes\mathfrak D, & i = 0  \\
e_i \otimes B,                       & i\ne 0
\end{cases}
\end{equation}
and consider induced grading on $H^2(L,L)$. $H_+^2(L,L)$ denotes a positive part of that
induced grading.

\begin{proposition}\label{5.2}
Let $L = S \otimes B + 1\otimes \mathfrak D, \> \mathfrak D \subseteq Der(B)$. Then:
\begin{enumerate}
\item if $S = W_1(n)$, then $H_+^2(L,L) \simeq H_+^2(S,S) \otimes B^{\mathfrak D} 
                            \>\oplus\> Der(B)^{\mathfrak D}$
\item if $S = sl(2)$, then $H_+^2(L,L) = 0$.
\end{enumerate}
\end{proposition}

\begin{proof}
(i) Proposition 1.1 ensures that Lemma 5.1 is applicable here if we put
$I = S \otimes B$ and $Q = 1\otimes \mathfrak D$. Using Theorem 4.4, (4.8) and (4.9)
and considering the action of $1\otimes \mathfrak D$ on appropriate cohomology groups
on the level of explicit cocycles, one gets
\begin{align}
E_2^{11} &\simeq H^1(S,S) \otimes H^1(\mathfrak D, B) \>\oplus\> 
                 H^1(\mathfrak D, Der(B)/\mathfrak D)                        \notag  \\
E_2^{02} &\simeq H^2(S,S) \otimes B^{\mathfrak D} \>\oplus\> Der(B)^{\mathfrak D} \>\oplus\>
                 Har^2(B,B)^{\mathfrak D} \>\oplus\> (Ker\,F)^{\mathfrak D}. \notag
\end{align}

The grading (5.2) induces a $\mathbb Z$-grading on each term of the spectral sequence 
(cf. \cite{F}).

The knowledge of $H^1(W_1(n),W_1(n))$, $H^2(W_1(n),K)$ (cf. \cite{B1} or \cite{D2}) 
and $H^2(W_1(n),W_1(n))$ (cf. \cite{DK}), allows to write down all 
nonzero graded components of $E_2 = E_2^{11} \oplus E_2^{02}$   
with respect to grading (5.2):
\begin{align}
& (E_2)_{-p^n}      \simeq
         H_{-p^n}^2 (W_1(n),W_1(n)) \otimes B^{\mathfrak D} \>\oplus\> 
         (Ker\,F)^{\mathfrak D}                                     \notag \\
& (E_2)_{-p^t}      \simeq
         H_{-p^t}^2 (W_1(n),W_1(n)) \otimes B^{\mathfrak D} \>\oplus\> 
         H_{-p^t}^1 (W_1(n),W_1(n)) \otimes H^1(\mathfrak D, B)     \notag \\       
& (E_2)_0           \simeq Der(B)^{\mathfrak D} \>\oplus\> Har^2(B,B)^{\mathfrak D} 
         \>\oplus\> 
         H^1(\mathfrak D, Der(B)/\mathfrak D)                       \notag \\
& (E_2)_{p^n - p^t} \simeq H_{p^n - p^t}^2 (W_1(n),W_1(n)) \otimes B^{\mathfrak D}  
                                                                    \notag \\
& (E_2)_{p^n}       \simeq Der(B)^{\mathfrak D}                     \notag
\end{align}
where $1 \le t \le n-1$.

The last two classes constitute $(E_2)_+$ and generated by cocycles 
$(1\otimes R_u) \circ \Theta_{\psi_t}$ where $u \in B^{\mathfrak D}$,
\begin{equation}\notag
\psi_t (e_i, e_j) = \begin{cases}
e_{p^n-2}, & i = -1, j = p^t -1   \\
0,         & \text{otherwise}
\end{cases}
\end{equation}
and $\phi_D$ where $D \in Der(B)^{\mathfrak D}$, respectively. 
(Strictly speaking, these cocycles are extended from the corresponding cocycles from 
$Z^2 (W_1(n) \otimes B, W_1(n) \otimes B)$ by letting them vanish on 
$W_1(n) \otimes B \wedge 1\otimes \mathfrak D$ and 
$1\otimes \mathfrak D \wedge 1\otimes\mathfrak D$).

According to Lemma 5.1, the corresponding $E_3$-term is
\begin{equation}\notag
(E_3^{02})_+ = Ker((E_2^{02})_+ \overset{(d_2^{02})_+}\longrightarrow (E_2^{21})_+).
\end{equation}

Theorem 4.4 shows that all $\mathfrak D$-invariant cohomology classes in 
$(E_2^{02})_+ \simeq H_+^2(S,S) \otimes B^{\mathfrak D} \oplus (Der(B))^{\mathfrak D}$
can be represented by $\mathfrak D$-invariant cocycles, what implies $(d_2^{02})_+ = 0$
and the positive part of the spectral sequence collapses in the relevant range.

(ii) Quite analogous (and simpler).
\end{proof}

\begin{remark}
In principle, one may compute the whole cohomology group 
$H^2 (S \otimes B + 1\otimes \mathfrak D, S \otimes B + 1\otimes \mathfrak D)$ 
by the following scheme: first, it is possible to evaluate $(Ker\,F)^{\mathfrak D}$
in the spirit of \S 2 or \S 3, and, particularly, to show that in this case the
$Q$-action commutes with $F$, what in its turn implies
\begin{equation}\notag
E_2^{02} = H^2(S \otimes B, S \otimes B)^{\mathfrak D} \oplus 
(Ker\,F \cap (H^2(S \otimes B) \otimes \mathfrak D)^{\mathfrak D}).
\end{equation}
Then the same reasoning as at the end of the proof of Proposition 5.2 shows that
$d_2^{02} = 0$ in general.
\end{remark}

As all cocycles constituting the basis of $H_+^2(L,L)$, being of the types
$(1\otimes R_u) \circ \Theta_{\psi_t}$ and $\phi_D$,
are possibly nonzero only on the ($-1$)\textit{st}, $1$\textit{st} and 
($p^t-1$)\textit{st} ($1 \le t \le n-1$) graded components of $L$, 
with values in the ($p^n-2$)\textit{th} graded component, each Massey product of two such
cocycles is obviously zero, and by Proposition 1.4 we have

\begin{theorem}\label{5.3}
Let $L$ be as in Proposition 5.2 with grading defined by (5.2). Let
$\mathfrak L$ be a filtered algebra whose associated graded algebra is isomorphic 
(as graded algebra) to $L$. Then
\begin{enumerate}
\item if $S = W_1(n)$, then $\mathfrak L$ is determined by brackets
\begin{equation}\notag
\curlybrack = \liebrack + \sum_{t=1}^{n-1} (1\otimes R_{u_t}) \circ \Theta_{\psi_t} + \phi_D
\end{equation}
for some $u_t \in B^{\mathfrak D}$ and $D \in Der(B)^{\mathfrak D}$
\item if $S = sl(2)$, then $\mathfrak L \simeq L$.
\end{enumerate}
\end{theorem}

Note that the basic cocycles in $H_+^2(L,L)$ are of the form $\Phi_d$ for appropriate
$d \in Der(O_1(n-1) \otimes B)$: $\phi_D = \Phi_{x^{p^{n-1}-1} \otimes D}$ and
$(1\otimes R_u) \circ \Theta_{\psi_t} = \Phi_{x^{p^{n-1}-1} \partial^{p^{t-1}} \otimes R_u},
1 \le t \le n-1$. Therefore the algebras appearing in part (i) of Theorem 5.3 are all of
the kind $\mathfrak L(A, D) + 1\otimes \mathfrak D$ 
for $A = O_1(n-1) \otimes B$ and $\mathfrak D \subseteq Der(A)$. 
(Note that from the Jacobi identity follows $[D,\mathfrak D] = 0$). 
Combining this fact (in the particular case $B = O_m$) with Theorem 1.5, we conclude:

\begin{corollary}\label{5.4}
(The ground field is algebraically closed of characteristic $p > 5$).

Let $\mathfrak L$ be a semisimple Lie algebra with a solvable maximal
subalgebra defining in $\mathfrak L$ a long filtration. Then either
$\mathfrak L \simeq sl(2) \otimes O_m + 1\otimes \mathfrak D$, or
$\mathfrak L \simeq \mathfrak L(O_m, D) + 1\otimes \mathfrak D$ for some
$m \in \mathbb N$, $D \in Der(O_m)$ and a solvable subalgebra $\mathfrak D$ in $Der(O_m)$
such that $[D,\mathfrak D] = 0$ and $O_m$ has no $\langle \mathfrak D, D \rangle$-invariant ideals.
\end{corollary}

\begin{remarks}
\hfill

(i)
The close inspection of Weisfeiler's results shows that if $\mathfrak L_0$ is a solvable maximal
subalgebra in Theorem 1.5, then after passing to the associated graded algebra,
$\mathfrak L_0$ goes to $\langle e_0, e_1 \rangle \otimes O_m + 1\otimes \mathfrak D$
(in the case $S = sl(2)$) or to $W_1(n)_0 \otimes O_m + 1\otimes \mathfrak D$ 
(in the case $S = W_1(n)$). Our computations of filtered deformations show that actually
$\mathfrak L_0$ coincides with these algebras (as they do not change under 
deformations).

(ii)
Since $[D, \mathfrak D] = 0$, the algebra $\langle \mathfrak D, D \rangle \subseteq Der(O_m)$
either coincides with $\mathfrak D$ (if $D\in \mathfrak D$), or is 
1-dimensional abelian extension of $\mathfrak D$ (if $D\notin \mathfrak D$).
\end{remarks}

So, to classify semisimple Lie algebras with a solvable maximal subalgebra occurring
in Theorem 1.5, it remains to describe algebras appearing in Corollary 5.4 up
to isomorphism and to identify them with the known semisimple Lie algebras. 
We accomplish this task in the next section.

\section{Classification of semisimple algebras $\mathfrak L(A, D) + 1\otimes \mathfrak D$}

The object of this section is the class of Lie algebras 
$\mathfrak L(A, D) + 1\otimes \mathfrak D$, where $1\otimes \mathfrak D$ acts on
$\mathfrak L(A, D)$ as on $W_1(1) \otimes A$ and $[D,\mathfrak D] = 0$.
The case $A = O_m$ is of particular importance.

All algebras throughout this section assumed to be finite-dimensional.

Note that consideration of dimensions immediately implies that no algebra of
the form $\mathfrak L(O_n, D)$ is isomorphic to some $sl(2) \otimes O_m$.

\begin{lemma}\label{6.1}
Each ideal of $\mathfrak L(A, D) + 1\otimes \mathfrak D$ is of the form
$\mathfrak L(I, D) + 1\otimes \mathfrak E$ where $I$ is a $\langle \mathfrak D, D \rangle$-invariant ideal
of $A$, $\mathfrak E$ is an ideal of $\mathfrak D$, and $\mathfrak E(A) \subseteq I$.
\end{lemma}

\begin{proof}
Let $\mathfrak I$ be an ideal of $\mathfrak L(A, D) + 1\otimes \mathfrak D$, 
then $\mathfrak I \cap \mathfrak L(A, D)$ is an ideal of $\mathfrak L(A, D)$.
Passing to the associated graded algebra (as in Proposition 1.3), 
we get that $gr(\mathfrak I \cap \mathfrak L(A, D))$
is an ideal of $W_1(1) \otimes A$. Either a direct calculation in $W_1(1)$,
or the general result of \cite{Ste}, yields that 
\begin{equation}\tag{6.1}
gr(\mathfrak I \cap \mathfrak L(A, D)) = W_1(1) \otimes I
\end{equation}
for some ideal $I$ of $A$.
Particularly, $e_{p-2}\otimes I \subset \mathfrak I \cap \mathfrak L(A,D)$.
Multiplying elements from $e_{p-2}\otimes I$ a necessary number of times
by $e_{-1}\otimes 1$, one gets $e_i\otimes I \subset \mathfrak I \cap \mathfrak L(A,D)$
for each $-1\le i \le p-2$, that is, $W_1(1)\otimes I \subseteq \mathfrak I \cap \mathfrak L(A,D)$.
Due to (6.1) this inclusion is actually an equality (of vector spaces): 
$\mathfrak I \cap \mathfrak L(A,D) = W_1(1)\otimes I$. 
Particularly, $W_1(1)\otimes I$ is closed under brackets $\curlybrack$
(cf. Definition in \S 1), what is equivalent to $D(I) \subseteq I$.

Now, taking an arbitrary element 
$\sum_{i=-1}^{p-2} e_i \otimes a_i + 1\otimes d \in \mathfrak I$, and multiplying it by
$e_0 \otimes 1$ and $e_{-1} \otimes 1$, we get 
$\sum_i ie_i \otimes a_i \in \mathfrak I \cap \mathfrak L(A,D)$ and 
$\sum_i e_{i-1} \otimes a_i \in \mathfrak I \cap \mathfrak L(A, D)$
respectively, showing therefore that all $a_i \in I$. This proves that
$\mathfrak I = \mathfrak L(I, D) + 1\otimes \mathfrak E$
for some subalgebra $\mathfrak E \subseteq \mathfrak D$. 
The rest of the conditions in the assertion follow immediately.
\end{proof}

\begin{lemma}\label{6.2}
Let a Lie algebra $\mathfrak L = \mathfrak L(A, D) + 1\otimes \mathfrak D$ 
is semisimple. Then the following hold:
\begin{enumerate}
\item $A \simeq \bigoplus_i O_{n_i}$ for some $n_i \in \mathbb N$, 
      each $O_{n_i}$ has no $\langle \mathfrak D, D \rangle$-invariant ideals
\item $\mathfrak L(A, D) \simeq \bigoplus_i S_i \otimes O_{m_i}$ for some $m_i \in \mathbb N$
      and simple Lie algebras $S_i$
\item $S_i \simeq \mathfrak L(O_{k_i}, d_i)$ for some $k_i \in \mathbb N$ and 
      $d_i \in Der(O_{k_i})$, $O_{k_i}$ has no $d_i$-invariant ideals.
\end{enumerate}
\end{lemma}

\begin{proof}
The proof merely consists of multiple applications of classical Block's results \cite{block-diff}.

By Lemma 6.1, $A$ has no $\langle \mathfrak D, D \rangle$-invariant nilpotent
ideals, i.e., $A$ is $\langle \mathfrak D, D \rangle$-semisimple in the terminology of Block 
\cite{block-diff}. 
According to \cite{block-diff}, Main Theorem and Corollary 8.3, $A$ is isomorphic to the direct sum
$\bigoplus_i O_{n_i}$ of reduced polynomial rings having no 
$\langle \mathfrak D, D \rangle $-invariant ideals. 
Hence $D = \sum_i D_i$, where each $D_i$ acts as derivation on $O_{n_i}$ and zero on 
$O_{n_j}, j \ne i$. Obviously
\begin{equation}\notag
\mathfrak L( \bigoplus_i O_{n_i}, \sum_i D_i ) \simeq \bigoplus_i \mathfrak L(O_{n_i}, D_i)
\end{equation}
and by Lemma 6.1 each minimal ideal of $\mathfrak L$ coincides with one of 
$\mathfrak L(O_{n_i}, D_i)$. Thus by \cite{block-diff}, Theorem 1.3, 
$\mathfrak L(O_{n_i}, D_i) \simeq S_i \otimes O_{m_i}$ for some simple Lie algebra $S_i$ and
$m_i \in \mathbb N$.

Applying Lemma 6.1 again, we see that each ideal of $\mathfrak L(O_{n_i}, D_i)$ is of the form
$\mathfrak L(I, D_i)$ for some $D_i$-invariant ideal $I$ of $O_{n_i}$, and by \cite{Ste}
each ideal of $S_i \otimes O_{m_i}$ is of the form $S_i \otimes J$ for some ideal $J$ of 
$O_{m_i}$. But $O_{m_i}^+$ is the greatest ideal of $O_{m_i}$
whence there is a greatest $D_i$-invariant ideal $I_i$ of $O_{n_i}$, 
$\mathfrak L(I_i, D_i) \simeq S \otimes O_{m_i}^+$, and
\begin{equation}\notag
\mathfrak L(O_{n_i}, D_i) / \mathfrak L(I_i,D_i) \simeq 
(S_i \otimes O_{m_i}) / (S_i \otimes O_{m_i}^+) \simeq S_i.
\end{equation}
It is easy to see that the left side here is isomorphic to
$\mathfrak L(O_{n_i}/I_i, d_i)$, \\ $d_i \in Der(O_{n_i}/I_i)$ being induced from $D_i$.
Since $S_i$ is simple, $O_{n_i}/I_i$ has no $d_i$-invariant ideals and again 
by Block's theorem, $O_{n_i}/I_i \simeq O_{k_i}$ for some $k_i \in \mathbb N$.
\end{proof}

Now we determine simple Lie algebras in the class $\mathfrak L(A, D)$.

\begin{lemma}\label{6.3}
(The ground field is perfect of characteristic $p > 3$).

$\mathfrak L = \mathfrak L(A, D)$ is simple if and only if $\mathfrak L \simeq W_1(n)$
for some $n \in \mathbb N$.
\end{lemma}

\begin{proof}
The ``if'' part contained in Proposition 1.2. So suppose that $\mathfrak L(A, D)$ is
simple. According to Lemmas 6.1 and 6.2, $A \simeq O_n$ for some $n \in \mathbb N$.
Hence $\mathfrak L$ has a subalgebra 
$\mathfrak L_0 = e_{-1} \otimes O_n^+ + \langle e_0, \dots, e_{p-2} \rangle \otimes O_n$
of codimension 1. Then by \cite{D1}, $\mathfrak L$ is isomorphic to either $sl(2)$
or $W_1(n)$, the first case is impossible by dimension consideration. 
\end{proof}

\begin{remarks} \hfill

(i) 
If the ground field is algebraically closed, one may deduce the assertion of the
Lemma from many other results in the literature, e.g. \cite{ree} 
(by utilizing the fact that algebras under consideration are Ree's algebras,
see remark after definition of $\mathfrak L(A, D)$ in \S 1), or \cite{K} or \cite{W}
(by noting that that $\mathfrak L_0$ is solvable).

(ii)
Combining Theorem 5.3(i) (with remark after it) and Lemma 6.3, we
recover the fact that each filtered deformation (with respect to the standard grading)
of $W_1(n)$ is isomorphic to $W_1(n)$. This fact is important in consideration of some
classes of Lie algebras with given properties of subalgebras or elements and was
proved by Benkart, Isaacs and Osborn in \cite{BIO}, \S 3 and
Dzhumadil'daev in \cite{D1}.
\end{remarks}

Now summarizing all our results, we obtain the final classification of the long
filtration case.

\begin{theorem}\label{6.4}
(The ground field is algebraically closed of characteristic $p > 5$).

$\mathfrak L$ is a semisimple Lie algebra with a solvable maximal
subalgebra defining in it a long filtration, if and only if either 
$\mathfrak L \simeq sl(2) \otimes O_m + 1 \otimes \mathfrak D$, or
$W_1(n) \otimes O_m \subset \mathfrak L \subset Der(W_1(n)) \otimes O_m + 1\otimes W_m$,
where $\mathfrak D$ in the first case, and $pr_{W_m}\mathfrak L$ in the second, 
are solvable subalgebras of $W_m$ such that $O_m$ 
has no $\mathfrak D$- or $pr_{W_m}\mathfrak L$-invariant ideals.
\end{theorem}

\begin{proof}

\textit{``only if'' part}. 
Summarizing results of Lemmas 6.2 and 6.3, we obtain that semisimple Lie
algebras of the form $\mathfrak L(A, D) + 1 \otimes \mathfrak D$ are exactly those 
whose socle is a direct sum of algebras $W_1(n) \otimes O_m$ for some $n, m \in \mathbb N$. 
By Corollary 5.4, these algebras (with solvable $\mathfrak D$), 
along with $sl(2)\otimes O_m + 1\otimes \mathfrak D$, exhaust all possible
semisimple Lie algebras with a solvable maximal subalgebra defining in it a long 
filtration. Obviously a socle of such algebra should consist of only one minimal
ideal, and the assertion desired follows.

\textit{``if'' part}. 
In the $sl(2)$ case it is evident that
$\mathfrak L_0 = \langle e_0, e_1 \rangle \otimes O_m + 1\otimes \mathfrak D$ 
is a solvable maximal subalgebra. 

In the $W_1(n)$ case, we have
$W_1(n) \otimes O_m \subset \mathfrak L \subset Der(W_1(n) \otimes O_m) \simeq
Der(W_1(n)) \\ \otimes O_m + 1 \otimes Der(O_m)$. 
By Proposition 1.2, identify $W_1(n)\otimes O_m$ with 
$\mathfrak L(O_1(n-1)\otimes O_m, \partial\otimes 1)$. By Theorem 3.1(ii),
$\mathfrak L = \mathfrak L(O_1(n-1)\otimes O_m, \partial\otimes 1) + 1\otimes\mathfrak D$
for some solvable subalgebra $\mathfrak D \subset Der(O_1(n-1)\otimes O_m)$ 
(a further elucidation
of the structure of $\mathfrak D$ is possible due to conditions imposed on $pr_{Der(O_m)}\mathfrak L$
and Theorem 4.1(iii), but we don't need it here).

Consider a maximal subalgebra $\mathfrak L_0$ containing a subalgebra
$\langle e_0, e_1, \dots, e_{p-2} \rangle \otimes O_1(n-1) \otimes O_m + 1\otimes \mathfrak D$.
Obviously 
$$
\mathfrak L_0 = e_{-1}\otimes I + \langle e_0, e_1, \dots, e_{p-2} \rangle \otimes O_1(n-1) \otimes O_m + 
1\otimes \mathfrak D
$$
for some $I \vartriangleleft O_1(n-1) \otimes O_m$. Since $O_1(n-1) \otimes O_m$ is a 
reduced polynomial ring itself, each its ideal is nilpotent, 
whence $\mathfrak L_0$ is solvable.
\end{proof}

\section*{Acknowledgements}

This paper was basically written at 1994 while being a Ph.D. student at Bar-Ilan University
under the guidance of Steve Shnider.

Back to those (and earlier) days, my thanks due to Askar Dzhumadil'daev 
for introducing me to the subject
and many useful advices, Steve Shnider for his generous support, attention, and significant 
conceptual improvements, Giora Dula for sharing his experience in spectral
sequences calculations, Malka Schaps for useful comments on a preliminary version
of the manuscript, and Baruch Muskat for help with the proof of Lemma 3.3.

I am also indebted to the referee for a very thorough reading of the manuscript
and numerous useful comments and corrections.


\begin{thebibliography}{BIO}

\bibitem[BIO]{BIO}
G.~M.~Benkart, I.~M.~Isaacs and J.~M.~Osborn,
\emph{Lie algebras with self-centralizing ad-nilpotent elements},
J.~Algebra \textbf{57} (1979), 279--309.

\bibitem[B1]{B1}
R.~E.~Block, \emph{On the extensions of Lie algebras},
Canad. J. Math. \textbf{20} (1968), 1439--1450.

\bibitem[B2]{block-diff}
R.~E.~Block,
\emph{Determination of the differentiably simple rings with a minimal ideal},
Ann. Math. \textbf{90} (1969), 433--459.

\bibitem[C]{C}
J.~L.~Cathelineau,
\emph{Homologie de degr\'e trois d'alg\`ebres de Lie simple 
d\'eploy\'ees \'etendues \`a une alg\`ebre commutative},
Enseign. Math. \textbf{33} (1987), 159--173.

\bibitem[D1]{D1}
A.~S.~Dzhumadil'daev, \emph{Simple Lie algebras with a subalgebra of codimension one},
Russian Math. Surv. \textbf{40} (1985), No. 1, 215--216.

\bibitem[D2]{D2}
A.~S.~Dzhumadil'daev,
\emph{Central extensions of the Zassenhaus algebra and their irreducible representations},
Math. USSR Sbornik \textbf{54} (1986), 457--474.

\bibitem[DK]{DK}
A.~S.~Dzhumadil'daev and A.~I.~Kostrikin, \emph{Deformations of the Lie algebra $W_1(m)$},
Proc. Steklov Inst. Math. \textbf{4} (1980), 143--158.

\bibitem[F]{F}
D.~B.~Fuchs, \emph{Cohomology of Infinite-Dimensional Lie Algebras},
Consultants Bureau, N.~Y., 1986.

\bibitem[GS]{GS}
M.~Gerstenhaber and S.~D.~Schack, \emph{Algebraic cohomology and deformation theory},
Deformations Theory of Algebras and Structures and Applications 
(M.~Hazewinkel and M.~Gerstenhaber, eds.), Kluwer, 1988, pp.~11--264.

\bibitem[Ha]{harrison}
D.~K.~Harrison, \emph{Commutative algebras and cohomology},
Trans. Amer. Math. Soc. \textbf{104} (1962), 191--204.

\bibitem[Ho]{holm}
T.~Holm, \emph{Hochschild cohomology rings of algebras $k[X]/(f)$},
Beitr\"age zur Algebra und Geometrie \textbf{41} (2000), No. 1, 291--301.

\bibitem[K]{K}
M.~I.~Kuznetsov,
\emph{Simple modular Lie algebras with a solvable maximal subalgebra},
Math. USSR Sbornik \textbf{30} (1976), 68--76.

\bibitem[L]{loday}
J.-L.~Loday, \emph{Cyclic Homology}, Springer, 1998.

\bibitem[PS]{PS}
A.~Premet and H.~Strade,
\emph{Simple Lie algebras of small characteristic. III. The toral rank $2$ case},
J.~Algebra \textbf{242} (2001), 236--337.

\bibitem[RS]{RS}
S.~Ramanan and R.~Sridharan, \emph{Resolutions for some filtered algebras},
Math. Ann. \textbf{148} (1962), 341--348.

\bibitem[Re]{ree}
R.~Ree, \emph{On generalized Witt algebras},
Trans. Amer. Math. Soc. \textbf{83} (1956), 510--546.

\bibitem[Ri]{riordan}
J.~Riordan, \emph{Combinatorial Identities}, Wiley, N.~Y., 1968.

\bibitem[Ste]{Ste}
I.~N.~Stewart, \emph{Tensorial extensions of central simple algebras},
J.~Algebra \textbf{25} (1973), 1--14.

\bibitem[Str]{Str}
H.~Strade, 
\emph{The classification of the simple modular Lie algebras: VI. Solving the final case},
Trans.~Amer.~Math.~Soc. \textbf{350} (1998), 2553--2628.

\bibitem[SW]{SW}
H.~Strade and R.~L.~Wilson,
\emph{Classification of simple Lie algebras over algebraically closed field of prime 
characteristic}, Bull. Amer. Math. Soc. \textbf{24} (1991), 357--362.

\bibitem[W]{W}
B.~Weisfeiler, \emph{On subalgebras of simple Lie algebras of characteristic $p>0$},
Trans. Amer. Math. Soc. \textbf{286} (1984), 471--503.

\bibitem[Z]{me-trans}
P.~Zusmanovich, \emph{Central extensions of current algebras},
Trans. Amer. Math. Soc. \textbf{334} (1992), 143--152.

\end{thebibliography}
\end{document}